\documentclass[11pt]{amsart}

\usepackage[utf8]{inputenc}
\usepackage[left=3cm, right=3cm, top=3cm, marginparwidth=2cm]{geometry}

\usepackage{mathtools}
\usepackage{amssymb}
\usepackage{bbold}

\usepackage[unicode,colorlinks=true,linkcolor=Brown,citecolor=Brown,urlcolor=Brown]{hyperref}
\usepackage[nameinlink,noabbrev]{cleveref}

\usepackage[dvipsnames]{xcolor}
\usepackage{tikz-cd}
\usepackage{enumerate}
\usepackage[textsize=tiny]{todonotes}

\usepackage[all]{xy}
\SelectTips{cm}{}
\newdir{ >}{{}*!/-10pt/\dir{>}}

\numberwithin{equation}{section} 

\swapnumbers 

\theoremstyle{plain}
\newtheorem{thm}[equation]{Theorem}
\newtheorem*{thm*}{Theorem}
\newtheorem{prop}[equation]{Proposition}
\newtheorem{lem}[equation]{Lemma}
\newtheorem{cor}[equation]{Corollary}

\theoremstyle{definition}
\newtheorem{dfn}[equation]{Definition}

\newtheorem{exas}[equation]{Examples}
\newtheorem{ntn}[equation]{Notation}

\theoremstyle{remark}
\newtheorem{hyp}[equation]{Hypothesis}
\newtheorem{cons}[equation]{Construction}
\newtheorem{rec}[equation]{Recollection}
\newtheorem{rem}[equation]{Remark}
\newtheorem*{AIstatement}{AI Statement}
\newtheorem*{ack}{Acknowledgments}

\AddToHook{env/thm/begin}{\crefalias{equation}{thm}}
\AddToHook{env/prop/begin}{\crefalias{equation}{prop}}
\AddToHook{env/lem/begin}{\crefalias{equation}{lem}}
\AddToHook{env/cor/begin}{\crefalias{equation}{cor}}

\AddToHook{env/dfn/begin}{\crefalias{equation}{dfn}}
\AddToHook{env/exa/begin}{\crefalias{equation}{exa}}
\AddToHook{env/ntn/begin}{\crefalias{equation}{ntn}}

\AddToHook{env/hyp/begin}{\crefalias{equation}{hyp}}
\AddToHook{env/cons/begin}{\crefalias{equation}{cons}}
\AddToHook{env/rec/begin}{\crefalias{equation}{rec}}
\AddToHook{env/rem/begin}{\crefalias{equation}{rem}}

\Crefname{thm}{Theorem}{Theorems}
\Crefname{prop}{Proposition}{Propositions}
\Crefname{lem}{Lemma}{Lemmas}
\Crefname{cor}{Corollary}{Corollaries}
\Crefname{dfn}{Definition}{Definitions}
\Crefname{exa}{Example}{Examples}
\Crefname{ntn}{Notation}{Notations}
\Crefname{hyp}{Hypothesis}{Hypotheses}
\Crefname{cons}{Construction}{Constructions}
\Crefname{rec}{Recollection}{Recollections}
\Crefname{rem}{Remark}{Remarks}

\newcommand{\nc}{\newcommand}
\nc{\dmo}{\DeclareMathOperator}

\dmo{\Ab}{Ab}
\dmo{\Add}{Add}
\dmo{\coev}{coev}
\dmo{\coker}{coker}
\dmo{\comp}{comp}
\dmo{\Compl}{Compl}
\dmo{\cname}{\rm c}
\dmo{\cone}{cone}
\dmo{\Der}{D}
\dmo{\dname}{\rm d}
\dmo{\End}{End}
\dmo{\Ext}{Ext}
\dmo{\fpname}{\rm{fp}}
\dmo{\hname}{\rm h}
\dmo{\Hom}{Hom}
\dmo{\Id}{Id}
\dmo{\im}{im}
\dmo{\incl}{incl}
\dmo{\Ker}{Ker}
\dmo{\kos}{kos}
\dmo{\Loc}{Loc}
\dmo{\Max}{Max_{\otimes}^{Serre}}
\dmo{\modname}{mod}
\dmo{\Mod}{Mod}
\dmo{\opname}{op}
\dmo{\res}{res}
\dmo{\SH}{SH}
\dmo{\Spc}{Spc}
\dmo{\Spch}{Spc^{{\rm h}}}
\dmo{\Spec}{Spec}
\dmo{\Spech}{Spec\homol}
\dmo{\supp}{supp}
\dmo{\Supp}{Supp}

\nc{\adj}{\dashv}
\nc{\AH}{\layer{\cA}{\cH}}
\nc{\AI}{\layer{\cA}{\cI}}
\nc{\bbN}{\mathbb{N}}
\nc{\bbQ}{\mathbb{Q}}
\nc{\bbZ}{\mathbb{Z}}
\nc{\Bmax}{\cB_{\max}}
\nc{\Bmin}{\cB_{\min}}
\nc{\cat}[1]{\mathcal{#1}}
\nc{\cA}{\cat{A}}
\nc{\cB}{\cat{B}}
\nc{\cC}{\cat{C}}
\nc{\cD}{\cat{D}}
\nc{\cF}{\cat{F}}
\nc{\cG}{\cat{G}}
\nc{\cH}{\cat{H}}
\nc{\cI}{\cat{I}}
\nc{\cJ}{\cat{J}}
\nc{\cK}{\cat{K}}
\nc{\cL}{\cat{L}}
\nc{\cM}{\cat{M}}
\nc{\cN}{\cat{N}}
\nc{\colim}{\mathop{\mathrm{colim}}}
\nc{\compact}{^{\cname}}
\nc{\compl}{^{\scriptscriptstyle\complement}}
\nc{\cP}{\cat{P}}
\nc{\cQ}{\cat{Q}}
\nc{\cS}{\cat{S}}
\nc{\cT}{\cat{T}}
\nc{\cU}{\cat{U}}
\nc{\CY}{\Compl_Y}
\nc{\Db}{\Der^{\mathrm{b}}}
\nc{\Dperf}{\Der^{\mathrm{perf}}}
\nc{\dual}{^{\dname}}
\nc{\ee}{\mathbb{e}}
\nc{\eg}{{\sl e.\,g.}}
\nc{\EI}{E_{\cI}}
\nc{\eY}{\ee_Y}
\nc{\ff}{\mathbb{f}}
\nc{\fp}{^{\fpname}}
\nc{\fY}{\ff_Y}
\nc{\gm}{\mathfrak{m}}
\nc{\gp}{\mathfrak{p}}
\nc{\gq}{\mathfrak{q}}
\nc{\Homcat}[1]{\Hom_{\cat #1}}
\nc{\homol}{^{\hname}}
\nc{\hook}{\hookrightarrow}
\nc{\hT}{\widehat{\cT}}
\nc{\ideal}[1]{\big\langle #1\big\rangle}
\nc{\ie}{{\sl i.e.}\ }
\nc{\Imax}{\cI_{\max}}
\nc{\into}{\mathop{\rightarrowtail}}
\nc{\inv}{^{-1}}
\nc{\isotoo}{\mathop{\,\overset{\sim}\longrightarrow}\,}
\nc{\isoto}{\mathop{\,\overset{\sim}\to}\,}
\nc{\layer}[2]{\bar{#1}_{#2}}
\nc{\LOC}{\mathsf{LOC}^\otimes}
\nc{\MhT}{\Modd\hT\dual}
\nc{\mhT}{\modd\hT\dual}
\nc{\Mid}{\,\big|\,}
\nc{\MK}{\Modd\cat K}
\nc{\mK}{\modd\cat K}
\nc{\ML}{\Modd\cat L}
\nc{\mL}{\modd\cat L}
\nc{\Modd}{\Mod\text{\!-}}
\nc{\modd}{\modname\text{\!-}}
\nc{\MT}{\Modd\cT\dual}
\nc{\mT}{\modd\cT\dual}
\nc{\nos}{\pi}
\nc{\onto}{\mathop{\twoheadrightarrow}}
\nc{\op}{^{\opname}}
\nc{\ouratop}[2]{\genfrac{}{}{0pt}{}{\scriptstyle #1}{\scriptstyle #2}}
\nc{\phih}{\varphi^{{\rm h}}}
\nc{\potimes}[1]{^{\otimes #1}}
\nc{\QB}{Q_{\cB}}
\nc{\QHI}{Q_{\cH,\cI}}
\nc{\QH}{Q_{\cH}}
\nc{\QI}{Q_{\cI}}
\nc{\QJ}{Q_{\cJ}}
\nc{\qquadtext}[1]{\qquad\textrm{#1}\qquad}
\nc{\quadtext}[1]{\quad\textrm{#1}\quad}
\nc{\RB}{R_{\cB}}
\nc{\restr}[1]{\mathord{\vert}_{#1}}
\nc{\RI}{R_{\cI}}
\nc{\SER}{\mathsf{SER}^\otimes}
\nc{\SET}[2]{\big\{ #1\Mid#2 \big\}}
\nc{\sminus}{\!\smallsetminus\!}
\nc{\SpchK}{\Spch(\cat K)}
\nc{\SpchL}{\Spch(\cat L)}
\nc{\SpchT}{\Spch(\cT\dual)}
\nc{\SpcK}{\Spc(\cat K)}
\nc{\SpcL}{\Spc(\cat L)}
\nc{\SpcT}{\Spc(\cat T\dual)}
\nc{\unit}{\mathbb{1}}
\nc{\WhI}{\Wh(\cI)}
\nc{\WhJ}{\Wh(\cJ)}
\nc{\Wh}{\W\homol}
\nc{\WI}{\W(\cI)}
\nc{\WJ}{\W(\cJ)}
\nc{\W}{W}
\nc{\xto}[1]{\xrightarrow{#1}}
\nc{\yon}{\mathrm{h}}

\title{Closed layers in tensor-triangular geometry}
\date{2026 November 11}

\author[P.\ Balmer]{Paul Balmer}
\address{\newline \indent UCLA Mathematics Department, Los Angeles, CA 90095-1555, USA \newline}
\email{balmer@math.ucla.edu, rkommerell@math.ucla.edu, emmyvr@math.ucla.edu}

\author[R.\ Kommerell]{Rhea Kommerell}

\author[E.\ Van Rooy]{Emmy Van Rooy}


\thanks{First-named author is supported by a travel grant of Simons Foundation International.}


\begin{document}

\begin{abstract}
In the absence of a theory of closed immersions in tensor-triangular geometry, we propose an alternative approach via quotients of module categories.
As an application, we prove invariance of homological residue fields under completion.
\end{abstract}

\maketitle

\vskip-\baselineskip
\vskip-\baselineskip

\section{Introduction}

Let us recall some commutative algebra. Let $I\subseteq R$ be a finitely generated ideal and~$\hat{R}=\lim_{n}R/I^n$ the completion of~$R$ along the closed subset~$Z=V(I)$ of~$\Spec(R)$. Then:
\begin{enumerate}[\indent(A)]
\item
\label{it:A}%
The ring homomorphism $R\to \hat{R}$ induces an isomorphism~$R/I\isoto \hat{R}/I\hat{R}$.
\end{enumerate}
As a consequence, the spectrum and the residue fields do not change on the locus of completion:
\begin{enumerate}[\indent(A)]
\setcounter{enumi}{1}
\item
\label{it:B}%
The induced map $\varphi\colon \Spec(\hat{R})\to \Spec(R)$ restricts to a bijection~$\varphi\inv(Z)\isoto Z$.
\smallbreak
\item
\label{it:C}%
For every point $\gp\in Z$, the homomorphism $R\to \hat{R}$ induces an isomorphism between the residue field of~$R$ at~$\gp$ and the residue field of~$\hat{R}$ at the unique point~$\hat\gp$ in~$\varphi\inv(\gp)$.
\end{enumerate}
The ring~$\hat R$ and statements~\eqref{it:B} and~\eqref{it:C} only depend on the closed subset~$Z$. On the other hand,~\eqref{it:A} is a somewhat deeper statement since $R/I$ depends on the choice of the ideal~$I$.

Sanders and the first-named author generalized~\eqref{it:B} to tensor-triangular geometry as follows (see \Cref{rec:completion}). Let $\cT$ be a `big' tensor triangulated category and~$\cT\dual=\cT\compact$ the subcategory of dualizable-compact objects. Let $Y\subseteq \SpcT$ be a Thomason subset of the tt-spectrum, \eg\ a closed subset with quasi-compact complement.
Consider the `$Y$\!-completion' functor~$\CY=(-)^{\wedge}_{Y}=[\eY,-]$ and the induced map on spectra
\begin{equation}
\label{eq:completion}%
\CY\colon \cT\dual\to \hT\dual
\qquadtext{and}
\varphi=\Spc(\CY)\colon \Spc(\hT\dual)\to \SpcT.
\end{equation}
Then $\varphi$ restricts to a bijection $\varphi\inv(Y)\isoto Y$. This is~\cite[Theorem~4.6]{BalmerSanders25}, and the same holds for homological spectra (\Cref{cor:Y-equiv-Spch}).
These results generalize~\eqref{it:B} to tt-geometry.

It is much harder to generalize~\eqref{it:A}, mainly because we lack a tensor-triangular analog of the quotient~$R/I$; see~\Cref{rem:not-T_Z}.
Case in point, the notion of `residue field' led to ruminations~\cite{BalmerKrauseStevenson19} and unnerving~\cite{Hyslopetal26pp,Coulembier26pp}. 
In this paper, we propose a workaround that allows us in particular to generalize~\eqref{it:C} to tensor-triangular geometry:
\begin{thm}[\Cref{thm:residue-completion}]
\label{thm:residue-completion-intro}%
Let $\cH\in \SpchT$ be a homological prime on~$Y$ and~$\hat\cH\in\Spch(\hT\dual)$ its unique preimage. Then the $Y$\!-completion functor~\eqref{eq:completion} induces an equivalence 
\begin{equation}
\label{eq:hres-equiv-intro}%
\frac{\MT}{\Loc(\cH)} \isotoo \frac{\MhT}{\Loc(\hat\cH)}
\end{equation}
of the corresponding homological residue fields.
\end{thm}

We discuss completion mostly as an illustration. Our higher ambition is to understand the correct tensor-triangular analog of closed subschemes.
As this remains elusive, we propose an indirect approach via what we call \emph{layers} (\Cref{dfn:layer}), \ie finite quotients
of the module category~$\cA=\MK$ over any rigid tensor triangulated category~$\cK$. More precisely, a layer 
\begin{equation}\label{eq:layer-intro}%
\AI\coloneq\frac{\MK}{\Loc{(\cI)}}
\end{equation}
is the Gabriel quotient of~$\cA$ by a Serre $\otimes$-ideal~$\cI\subseteq\mK$ of finitely presented $\cK$-modules.

Although loosely inspired by the quotients~$R/I$ of commutative algebra, the above notion of a layer is actually much broader.
For instance, homological residue fields~\eqref{eq:hres-equiv-intro} are layers -- the \emph{minimal} ones.
More strikingly, for every quasi-compact \emph{open} subset~$U\subseteq \SpcK$, the module category~$\Modd{\cK(U)}$ of the open piece~$\cK(U)$ is also a layer as in~\eqref{eq:layer-intro}; see~\cite[Proposition~2.9]{Balmer20a}. To understand what the \emph{closed} layers are, let us discuss some spectra.

\begin{hyp}
\label{hyp-steel}%
To keep things simple, assume that the homological spectrum and the tt-spectrum coincide:~$\SpchK=\SpcK$. We drop this assumption after the Introduction.
\end{hyp}

Every layer~$\AI$ as in~\eqref{eq:layer-intro} has a natural `spectrum' (\Cref{dfn:spec-layer}) that identifies with the following subset of the spectrum of~$\cK$:
\begin{equation}
\label{eq:spc(AI)-intro}%
\WI=\SET{\cH\in\SpchK}{\cI\subseteq \cH}.
\end{equation}
The notation~$W$ is meant to evoke the closed subset~$V(I)=\SET{\gp\in\Spec(R)}{I\subseteq \gp}$ of commutative algebra while reminding us that we are in a broader setting: $\WI$ need not be a closed subset in general. Here are the actual possibilities.
\begin{thm}[{\Cref{sec:procon}}]
\label{thm:procon-intro}%
The spectrum~$\WI$ is a proconstructible subset (\Cref{rec:constructible}) of~$\SpcK$. Conversely, every proconstructible subset is the spectrum of some layer.
\end{thm}
Proconstructible subsets include singletons and quasi-compact open subsets, as illustrated above by the residual layer~$\AH$ (whose spectrum is just~$\{\cH\}$) and the open layer~$\Modd{\cK(U)}$ (whose spectrum is~$U$).
Closed layers are simply those layers whose spectrum is closed.

In commutative algebra, several ideals~$I\subseteq R$ define the same closed subset~$Z=V(I)$. Among them, the radical~$\sqrt{I}=\cap_{\gp\in Z}\,\gp$ is the largest.
Similarly, many of our Serre $\otimes$-ideals~$\cI$ define the same proconstructible subset~$Y=\WI$ of~$\SpcK$. 
Among them, the maximal one is given by $\Imax(Y)=\cap_{\cH\in Y}\,\cH$; see~\Cref{prop:Imax}.

Remarkably, there is also a `lower bound' when $Y$ is a Thomason subset (see~\Cref{sec:Thom-layers}). This fact does not seem to have an analog in commutative algebra. It is expressed in terms of `large' subcategories of~$\MK$.
Define~$\Bmin(Y)$ to be the kernel of the functor
\[
 \cA=\MK\onto \Modd{\cK_Y}
\]
given by restriction from~$\cK$ to the tt-ideal~$\cK_Y=\SET{c\in\cK}{\supp(c)\subseteq Y}$ of objects supported on~$Y$. This localizing subcategory $\Bmin(Y)$ is typically \emph{not} of finite type. (See~\Cref{prop:Bmin-not-finite-type}.) We show that $\Bmin(Y)$ is a $\otimes$-ideal, which implies that $\Modd{\cK_Y}$ inherits a tensor product.

\begin{thm}[{\Cref{thm:Bmin}}]
\label{thm:bookend-intro}%
Let $Y\subseteq\SpcK$ be a Thomason subset and let $\Bmin(Y)$ be the above localizing $\otimes$-ideal of~$\cA$.
Let~$\cI$ be a Serre $\otimes$-ideal of~$\cA\fp=\mK$. Then the spectrum $\WI$ of the layer~$\bar\cA_\cI$ is contained in~$Y$ if and only if
\[
\Bmin(Y)\subseteq \Loc(\cI).
\]
In that case, the layer $\AI$ is a Gabriel quotient of~$\Modd{\cK_Y}$ by a localizing $\otimes$-ideal.
\end{thm}

It follows that the layers that live `over~$Y$' only depend on the tt-ideal~$\cK_Y$. To see this, consider $F\colon \cK\to \cL$ a tt-functor and~$\varphi=\Spc(F)\colon \SpcL\to \SpcK$ the associated map. For a Thomason subset~$Y\subseteq \SpcK$, we say that $F$ is an \emph{equivalence over~$Y$} if~$F$ restricts to an equivalence $\cK_Y\isoto \cL_{\varphi\inv(Y)}$.
For example, the $Y$-completion functor~$\cT\dual\to \hT\dual$ of~\eqref{eq:completion} is an equivalence over~$Y$, by~\cite[Theorem~4.6]{BalmerSanders25}.
We can now state our generalization of~\eqref{it:A} to tensor-triangular geometry; this occupies our central \Cref{sec:Y-invariance}.
\begin{thm}[{\Cref{thm:main}}]
\label{thm:main-intro}%
Let $F\colon \cK\to \cL$ be a tt-functor that is an equivalence over a Thomason subset~$Y$. Then $F$ induces a bijection between the Serre $\otimes$-ideals~$\cI$ of~$\mK$ whose spectrum~$\WI$ is contained in~$Y$ and the Serre $\otimes$-ideals~$\cJ$ of~$\mL$ whose spectrum~$\WJ$ is contained in~$\varphi\inv(Y)$. Moreover if~$\cI$ corresponds to~$\cJ$ under this bijection, then~$F$ yields an equivalence of the associated layers
\begin{equation*}
\frac{\MK}{\Loc(\cI)}\isotoo\frac{\ML}{\Loc(\cJ)}.
\end{equation*}
\end{thm}

The proof of this result relies on the following characterization of finite type $\otimes$-localizations of~$\MK$, that could be of independent interest.
\begin{thm}[{\Cref{thm:finite-type}}]
\label{thm:finite-type-intro}%
Let $\cB$ be a localizing~$\otimes$-ideal of~$\cA=\MK$. Then $\cB$ is generated by~$\cB\cap \cA\fp$ if and only if the $\otimes$-unit~$\bar{\unit}$ in~$\cA/\cB$ is finitely presented.
\end{thm}

As an illustration of our methods, we show in~\Cref{sec:regularity} how residual regularity at points of~$Y$ (in the sense of~\cite{VanRooy26}) is preserved along a tt-functor that is an equivalence over~$Y$.
As indicated above, this applies in particular to~$Y$-completion~$\cT\dual\to \hT\dual$.
\begin{thm}[{\Cref{thm:regularity}}]
\label{thm:regularity-intro}%
Let $F\colon \cK\to \cL$ be a tt-functor that is an equivalence over a Thomason subset~$Y$. Let $\varphi=\Spc(F)$. Let~$\cQ \in \SpcL$ such that $\cP\coloneq\varphi(\cQ)$ belongs to~$Y$. Then $\cL$ is residually regular at $\cQ$ if and only if $\cK$ is residually regular at~$\cP$. 
\end{thm}

In the final \Cref{sec:commutative-rings}, we consider the example of~$\cK=\Dperf(R)$ for a commutative ring~$R$. We show that the set of ideals of~$R$ injects into the set of Serre $\otimes$-ideals of $\mK$ in a way that is compatible with vanishing loci, further justifying the notation~$\WI$. This injection fails to be surjective already for the ring~$R=\bbZ$ of integers.

\begin{AIstatement}
The paper has been written entirely by the authors.
All new ideas and all new proofs are their own.
\end{AIstatement}

\begin{ack}
The authors are grateful to Jordan Williamson for references, and thank Beren Sanders and Greg Stevenson for comments.
\end{ack}


\tableofcontents

\section{Preliminaries and notation}
%

%
\begin{hyp}
\label{hyp:K-T}%
Throughout the paper,~$\cK$ denotes an essentially small tensor triangulated category that we assume idempotent-complete and rigid. Typically, $\cK$ is the subcategory of compact objects (or equivalently, dualizable objects) in a `big' \emph{rigidly-compactly generated} tensor triangulated category~$\cT$ in the sense of~\cite{BalmerFavi11}:
\begin{equation}
\label{eq:K-T}%
\cK=\cT\compact=\cT\dual.
\end{equation}
Such `big' tt-categories~$\cT$ are extremely common in topology, algebra and geometry. They are called `unital algebraic stable homotopy categories' in~\cite{HoveyPalmieriStrickland97}, where the reader can find a multitude of examples.
From the perspective of applications, there is no restriction in assuming the existence of a `big'~$\cT$ satisfying~\eqref{eq:K-T}, as $\cK$ can be Ind-completed in the presence of a model. We use~$\cK$ to emphasize that our theory only depends on small objects.
\end{hyp}

\begin{rec}
\label{rec:dual}%
Rigidity of $\cK$ implies that every object $c \in \cK$ admits a dual~$c^{\vee} \in \cK$ with an adjunction~$(c \otimes -) \adj ( c^{\vee}\otimes -)$. Its unit is denoted~$\coev_c\colon \unit \to c^{\vee} \otimes c$. As $c \otimes \coev_c \colon c \to c \otimes c^{\vee} \otimes c$ is a split monomorphism, every object~$c$ is a direct summand of~$c\otimes c^\vee\otimes c$.
\end{rec}

\begin{rec}
\label{rec:Yoneda}%
See~\cite[\S\,2]{BalmerKrauseStevenson19}.
We denote the category of right $\cK$-modules by
\[
\cA=\MK=\Add(\cK\op,\Ab).
\]
It is a locally coherent Grothendieck abelian category. Its abelian subcategory of finitely presented objects~$\cA\fp=\mK$, also known as the \emph{Freyd envelope} of~$\cK$ (\cite[Chapter~5]{Neeman01}), is Frobenius and its projective-injectives are the representables $\yon(c)$ for $c\in \cK$, where
\begin{equation*}
    \yon\colon \cK \hook \mK\hook \MK=\cA.
\end{equation*}
is the Yoneda embedding.
By Day convolution, $\cA$ inherits a tensor, cocontinuous in each variable and making $\yon$ a tensor functor. The category~$\cA$ has enough injectives and
\begin{equation}
\label{eq:inj=>flat}%
\textit{every injective in~$\cA$ is $\otimes$-flat.}
\end{equation}
See~\cite[\S2]{BirdWilliamson26}, or~\cite[Corollary~1.9]{Krause00} and~\cite[Proposition~A.14]{BalmerKrauseStevenson19} when~\eqref{eq:K-T} holds.
\end{rec}

\begin{rec}
\label{rec:F!-F*}%
Given $\cK$ and $\cL$ two essentially small triangulated categories, any exact functor $F \colon \cK \to \cL$ induces a functor $F_!\colon \MK \to \ML$, given by left Kan extension
\[
F_!(M) = \colim_{\substack{c\in\cK,\  \yon(c) \to M}} \yon(F(c)).
\]
In particular, $F_!\circ\yon = \yon\circ F$. This functor has a right adjoint $F^\ast \colon \ML \to \MK$, given by restriction along~$F$. Both $F_!$ and $F^{\ast}$ are exact cocontinuous functors and $F_!$ preserves finitely presented objects. See~\cite[\S 1]{Krause05}. When $\cK$ and $\cL$ are moreover tt-categories as in \Cref{hyp:K-T} and $F$ is a tt-functor, then $F_!$ is symmetric monoidal as well, and the adjunction $F_! \adj F^{\ast}$ satisfies a projection formula: $F^{\ast}(N \otimes F_!(M)) \cong F^{\ast}(N) \otimes M$. 
\end{rec}

\begin{rec}\label{rec:Spch}%
The \emph{homological spectrum of~$\cK$} is defined as $\SpchK=\Max(\cA\fp)$, where $\Max$ denotes the collection of maximal Serre $\otimes$-ideals. 
The latter are called \emph{homological primes}. See~\cite{Balmer20a}.
There is a well-defined surjection
\begin{equation*}
    \nos\colon \SpchK \onto \SpcK
\end{equation*}
sending a homological prime~$\cH\subset\mK$ to $\nos(\cH)=\yon\inv(\cH)\subset\cK$.
The surjection~$\nos$ is very often bijective, see~\cite[\S\,5]{Balmer20a}. However, $\nos$ can fail to be injective, as observed by~\cite{Hyslopetal26pp,Coulembier26pp}.
One sometimes says that $\cK$ has the \emph{steel property} if~$\nos$ is a bijection.
\end{rec}

\begin{rec}
\label{rec:A/I}%
See~\cite[\S\,2]{BalmerKrauseStevenson19}.
For every Serre $\otimes$-ideal~$\cI\subseteq \cA\fp=\mK$ we consider the Gabriel quotient by the localizing subcategory generated by~$\cI$ (which is always $\otimes$-ideal)
\begin{equation} 
\label{eq:AI}%
\begin{tikzcd} 
	{\cA = \MK} \\
	{\AI\coloneq\frac{\MK}{\Loc(\cI)}}
	\arrow[""{name=0, anchor=center, inner sep=0}, "{Q_\cI}"', shift right=3, two heads, from=1-1, to=2-1]
	\arrow[""{name=1, anchor=center, inner sep=0}, "{R_\cI}"', shift right=3, hook, from=2-1, to=1-1]
	\arrow["\adj"{anchor=center}, draw=none, from=0, to=1]
\end{tikzcd}
\end{equation}
The Grothendieck category~$\AI$ remains locally coherent and its subcategory of finitely presented objects~$\AI\fp$ identifies with~$(\mK)/\cI$, the quotient of~$\cA\fp$ by the chosen~$\cI$. In~\eqref{eq:AI}, the functor~$\RI$ is the fully faithful right adjoint to the exact $\otimes$-functor~$\QI$. The adjunction $\QI \adj \RI$ satisfies a projection formula:~$\RI(N \otimes \QI(M)) \cong \RI(N) \otimes M$. 
We often write $\bar{x}$ instead of~$\QI(x)$ for $x\in\cA$, or instead of~$\QI(\yon(x))$ for~$x\in \cK$. Generalizing~\eqref{eq:AI}, given two Serre $\otimes$-ideals $\cI\subseteq\cH$, we have an intermediate $\QHI\colon \AI\onto \AH$ such that $\QH\cong \QHI\circ\QI$.
\end{rec}

\begin{dfn}
\label{dfn:layer}%
We call the above finite Gabriel quotient~$\AI=\cA/\Loc(\cI)$ the \emph{layer of $\cK$ with respect to $\cI$}, or the \emph{$\cI$-layer of~$\cK$} for short.
\end{dfn}

\begin{rec}
\label{rec:E_I}%
The injective hull of the unit~$\bar{\unit}\into \bar{E}_{\cI}$ in~$\AI$ yields an injective object~$\EI\coloneq R_{\cI}(\bar{E}_{\cI})$ in~$\cA=\MK$ and a morphism~$\eta\colon \unit\to \EI$. Since $\QI\circ R_{\cI}=\Id$ we have $\QI(\EI)=\bar{E}_{\cI}$, justifying the notation.
(When~\eqref{eq:K-T} holds, this injective~$\EI\in\cA$ is the image of a pure-injective in~$\cT$ also denoted~$\EI$ in~\cite{BalmerKrauseStevenson19}.)
Thanks to~\eqref{eq:inj=>flat}, the object~$\EI$ is $\otimes$-flat in~$\cA$ and it follows that it is a weak ring ($\EI\otimes\eta\colon \EI\to \EI\otimes\EI$ has a retraction). We have 
\begin{equation}
\label{eq:EI}%
\Loc(\cI)=\Ker(\EI\otimes-\colon \cA\to \cA) 
\quadtext{and}
\cI=\Ker(\EI\otimes-)\cap \cA\fp.
\end{equation}
This is proved under~\eqref{eq:K-T} in~\cite{BalmerKrauseStevenson19} but the proof works in general; see~\cite[\S\,2]{BirdWilliamson26}.
\end{rec}

\begin{ntn}
We follow~\cite{VanRooy26} and depart from earlier literature by writing
\[
\kappa(\cH)
\]
for the injective object~$\EI$ when~$\cI=\cH$ is a homological prime~$\cH\in\SpchK$.
This is justified by~\cite{BalmerCameron21} where these objects are computed to be usual residue fields $\kappa(x)$ in algebraic geometry (not their injective hulls!) and to be Morava $K$-theories $K(n)$ in topology.
\end{ntn}

We believe that the following result is new. The `if' direction is somewhat surprising: One can decide whether a localizing $\otimes$-ideal~$\cB$ is of finite type solely in terms of~$\cA/\cB$.

\begin{thm}
\label{thm:finite-type}%
Let $\cB$ be a localizing $\otimes$-ideal subcategory of~$\cA=\MK$. Then $\cB$ is of finite type, \ie $\cB=\Loc(\cB \cap \cA\fp)$, if and only if the $\otimes$-unit~$\bar\unit$ in~$\cA/\cB$ is finitely presented.
\end{thm}
\begin{proof}
The forward direction is classical~\cite[Theorem~2.6]{Krause97}: $\cA/\cB$ remains locally coherent and~$\QB\colon \cA\onto \cA/\cB$ preserves finitely presented objects, including here~$\unit_{\cA}=\yon(\unit_{\cK})\in\mK=\cA\fp$.
The reverse direction is the interesting one. It suffices to prove that $\QB$ preserves finitely presented objects. Indeed, the usual argument using the adjunction $\QB\adj \RB$ and the fact that $\cA$ is locally finitely presented shows that $\RB\colon \cA/\cB\to \cA$ preserves filtered colimits, and this is sufficient to guarantee that $\cB$ is of finite type by~\cite[Lemma~2.4]{Krause97}. The little miracle here is that $\QB$, being a tensor-functor, preserves dualizable objects. It is easy to show in general that if~$\bar\unit$ is finitely presented in~$\cA/\cB$ then any dualizable~$d$ is finitely presented as well, \ie $\Hom(d,-)$ preserves filtered colimits: use $\Hom(d,-)\cong\Hom(\bar\unit,d^\vee\otimes-)$ and  the fact that~$d^\vee\otimes-$ commutes with all colimits. Now in~$\cA$ every finitely presented object is a cokernel of a map between dualizable ones, and~$\QB$ is exact. The result follows easily since a cokernel of a morphism of finitely presented objects remains finitely presented.
\end{proof}

\begin{rem}
\label{rem:not-T_Z}%
The concept of `closed immersion' has not been extended to tensor-triangular geometry, in the sense that there is no purely tensor-triangular construction that would produce the derived category~$\Der(Z)$ of a closed subscheme~$Z\hook X$ when applied to~$\cT=\Der(X)$. In particular, the naive guess~$\cT_Z$ produces~$\Der_Z(X)$, which only depends on the closed \emph{subset}~$Z$, as opposed to~$\Der(Z)$ that depends on the sheaf of ideals defining~$Z$. In fact, no closed subscheme of $X$ can have $\Der_Z(X)$ as its derived category because $\Der_Z(X)$ is not rigidly-compactly generated except in trivial cases. See a detailed discussion in~\cite[Remark~2.12]{BalmerSanders25}.
\end{rem}

\begin{rec}
\label{rec:eY-fY}%
Let $Y \subseteq \SpcK$ be a Thomason subset and consider the tt-ideal $\cK_Y = \{ a \in \cK \mid \supp(a) \subseteq Y\}$. 
When~\eqref{eq:K-T} holds, we can construct in~$\cT$ the idempotents~$\eY$ and~$\fY$ of~\cite{BalmerFavi11}. We explain how their images in $\cA=\MK$ can be constructed in the absence of~$\cT$.
\Cref{rec:F!-F*} for the inclusion $i \colon \cK_{Y} \hook \cK$ and the quotient $q \colon \cK \onto \cK/\cK_Y$ gives
\[\begin{tikzcd}[row sep=1.5em]
	{\MK_Y} && \MK \\
	\MK && {\Modd{\cK/\cK_Y}}
	\arrow[""{name=0, anchor=center, inner sep=0}, "{i_!}"', shift right=3, from=1-1, to=2-1]
	\arrow[""{name=1, anchor=center, inner sep=0}, "{q_!}"', shift right=3, from=1-3, to=2-3]
	\arrow[""{name=2, anchor=center, inner sep=0}, "{i^*}"', shift right=3, from=2-1, to=1-1]
	\arrow[""{name=3, anchor=center, inner sep=0}, "{q^*}"', shift right=3, from=2-3, to=1-3]
	\arrow["\adj"{anchor=center}, draw=none, from=0, to=2]
	\arrow["\adj"{anchor=center}, draw=none, from=1, to=3]
\end{tikzcd}\]
Using these adjunctions, we may define the following objects in $\MK$:
\begin{equation*}
    \label{equ:eY-fY}
    \eY \coloneq i_! i^* (\unit) \qquadtext{and} \fY \coloneq q^* q_!(\unit). 
\end{equation*}
The explicit formulas for these adjunctions give
\begin{equation}\label{eq:eY-fY-colims}
    \eY = \colim_{\substack{a \to \unit \textrm{ in }\cK\\[.1em] \textrm{ with }a \in \cK_Y}} \yon(a)
\qquadtext{and}
\fY = \Hom_{\cK/\cK_Y}(q(-), \unit) =\!\! \colim_{\substack{\sigma \colon\unit \to b \textrm{ in }\cK\\[.1em] \textrm{ with }\cone(\sigma) \in \cK_Y}} \yon(b).
\end{equation}
Since $\eY$ and $\fY$ are filtered colimits of $\otimes$-flat objects, they are themselves $\otimes$-flat. 
By \cite[Proposition 2.10]{BIK15}, the counit $\eY \to \unit$ and unit $\unit \to \fY$ fit into a long exact sequence in~$\cA$
\begin{equation}\label{eq:LES-eY-fY}
   \dots \to  \Sigma^{-1} \fY \to \eY \to \unit \to \fY \to \Sigma \eY \to \dots 
\end{equation}
Using the projection formula for the adjunction $q_! \adj q^*$ (\Cref{rec:F!-F*}), we see that
\begin{equation}
\label{eq:fY-tens-KY=0}%
\fY\otimes \yon(a)\cong q^*(q_!(\yon(a)))\cong q^*(\yon(q(a)))=0\qquad\textrm{ for all }a\in\cK_Y.
\end{equation}
This gives $\eY \otimes \fY = 0$ by~\eqref{eq:eY-fY-colims}. Hence $\eY \otimes \eY \isoto \eY$ {and}~$\fY \isoto \fY \otimes \fY$
by~\eqref{eq:LES-eY-fY} and $\otimes$-flatness of $\eY$ and $\fY$. So $\eY$ is a $\otimes$-idempotent co-ring and $\fY$ a $\otimes$-idempotent ring in~$\cA$.
\end{rec}
\begin{lem}[{Compare~\cite[Lemma 3.8]{BarthelHeardSanders23a}}]
\label{lem:eY-fY-H}%
With notation as in \Cref{rec:eY-fY}.
\begin{enumerate}[\rm(a)]
\item
\label{it:eY-fY-H-in}%
 We have $\SET{\cH\in\SpchK}{\QH(\fY)=0}=\SET{\cH\in\SpchK}{\QH(\eY)=\bar\unit}=\nos\inv(Y)$.
\smallbreak
\item 
\label{it:eY-fY-H-out}%
We have $\SET{\cH\in\SpchK}{\QH(\fY)=\bar\unit}=\SET{\cH\in\SpchK}{\QH(\eY)=0}=\nos\inv(Y)\compl$.
\end{enumerate}
\end{lem}
\begin{proof}
There are exactly two cases: If $\cH\in\nos \inv(Y)$ then there exists an object $a \in \cK_Y$ with $\nos(\cH) \in \supp(a)$, meaning $\bar{a} \neq 0$ in~$\AH$. But $\overline{\fY}\otimes\bar{a}=0$ in~$\AH$ because of~\eqref{eq:fY-tens-KY=0}, forcing~$\overline{\fY}=0$. This gives~$\overline{\eY}\cong\bar\unit$ by~\eqref{eq:LES-eY-fY}. The other case is $\cH\in\nos\inv(Y)\compl$. Then $\nos(\cH)\notin\supp(a)$ for all~$a\in\cK_Y$. Hence $\bar{a}=0$ for all~$a\in \cK_Y$ and it follows from~\eqref{eq:eY-fY-colims} that $\overline{\eY}=0$ as well. Again, \eqref{eq:LES-eY-fY} guarantees that $\overline{\fY}\cong\bar\unit$ in this case. This gives the result.
\end{proof}

\section{Spectra of layers}
\label{sec:procon}%

We keep~$\cK$ as in~\Cref{hyp:K-T} and extend the definition of the homological spectrum to every layer. We characterize the corresponding subsets of~$\SpcK$ topologically.

\begin{dfn}
\label{dfn:spec-layer}%
By extension of \Cref{rec:Spch}, the \emph{homological spectrum of the layer~$\AI$} (see~\Cref{dfn:layer}) refers to the set $\Max(\AI\fp)$ of maximal Serre $\otimes$-ideals of~$\AI\fp$. Since $\cA\fp\onto\AI\fp$ is the quotient by the chosen Serre $\otimes$-ideal~$\cI$, we have an injection $\Max(\AI\fp)\hook\Max(\cA\fp)=\SpchK$ given by~$\bar\cH\mapsto \QI\inv(\bar\cH)$. This identifies~$\Max(\AI\fp)$ with the following subset of~$\SpchK$, that we also call the \emph{homological spectrum of~$\AI$}:
\[
\WhI\coloneq\SET{\cH\in\SpchK}{\cI\subseteq\cH}.
\]
Its image under~$\nos\colon \SpchK\onto \SpcK$ is \emph{the spectrum of~$\AI$}
\[
\WI\coloneq\nos(\WhI).
\]
So the subset~$\WI\subseteq\SpcK$ consists of the primes~$\cP\in\SpcK$ where \emph{there exists} a homological prime~$\cH$ above~$\cP$ such that $\cI\subseteq\cH$, but there is no guarantee that this holds for all~$\cH\in\nos\inv(\cP)$.
When $\cK$ has the steel property, this distinction is irrelevant.
\end{dfn}

\begin{lem}
\label{lem:H-in-Spch(A/I)}%
Let $\cI\subseteq \mK$ be a Serre $\otimes$-ideal and let~$\cH\in\SpchK$. Then $\cH$ belongs to the homological spectrum $\WhI$ of the~$\cI$-layer~$\AI$ if and only if~$\kappa(\cH)\otimes\EI\neq 0$ in~$\MK$, that is, by~\eqref{eq:EI}, if and only if $\kappa(\cH)$ does not belong to~$\Loc(\cI)$.
\end{lem}
\begin{proof}
If $\cH$ belongs to $\WhI$, then $\cI \subseteq \cH$ and the Gabriel quotient $\QH$ factors via~$\AI$ as $\QH\cong \QHI\circ\QI$. Since $\QH(\kappa(\cH))=\bar{E}_{\cH}\neq 0$ this forces~$\QI(\kappa(\cH))\neq 0$, \ie $\kappa(\cH)\notin\Loc(\cI)$.

Conversely, suppose that $\kappa(\cH) \otimes \EI \neq 0$. Then, $\cH'\coloneq\Ker(\kappa(\cH) \otimes \EI\otimes -) \cap \cA\fp$ is a proper Serre $\otimes$-ideal of $\mK$ containing both $\cI$ and $\cH$, which is maximal. Hence $\cH=\cH'\supseteq\cI$. 
\end{proof}

\begin{rec}
\label{rec:nilpotence}%
Let $f\colon \unit\to X$ be a morphism in~$\cA$ with~$X$ $\otimes$-flat. 
The Nilpotence Theorem~\cite[Corollary 4.2]{Balmer20b} says that if $\bar f=0$ in~$\AH$ for all~$\cH\in \SpchK$ then $f\potimes{n}=0$ for~$n\gg1$.
If $f$ goes to zero in~$\AH$ only for those~$\cH\in \WhI$ then $f$ is $\otimes$-nilpotent in~$\AI$. One can either adapt the proof in \textsl{loc.\,cit.} to~$\AI$, or apply the original Nilpotence Theorem to $\eta\otimes f\colon \unit\to \EI\otimes X$ in~$\cA$, which vanishes also in~$\AH$ for~$\cH\notin\WhI$ by~\Cref{lem:H-in-Spch(A/I)}.
\end{rec}

\begin{rem}
\label{rem:piV}%
The spectrum $\SpcK$ is a spectral topological space. The homological spectrum $\SpchK$ is only a \emph{set}. Of course, one can pull back the topology of the former to the latter via~$\nos\colon \SpchK\onto \SpcK$ but by~\cite[Proposition~4.5]{BarthelHeardSanders23a} the resulting space $\SpchK$ is spectral if and only if~$\nos$ is a bijection. Hence we prefer to state topological properties of~$\WhI$ in terms of its image~$\WI$ in~$\SpcK$, particularly for properties involving spectral spaces. 
\end{rem}

\begin{rec}
\label{rec:constructible}%
In a spectral space, the constructible subsets generate a topology, called the \emph{constructible topology} or \emph{patch topology}. 
In the case of~$\SpcK$, a constructible subset is a finite union of sets of the form~$C(a,b)\coloneq\supp(a)\compl\cap \supp(b)=\SET{\cP\in\SpcK}{a\in\cP\textrm{ and }b\notin\cP}$ for~$a,b\in\cK$. Since $C(a,b)\cap C(a',b')=C(a\oplus a',b\otimes b')$ the~$\{C(a,b)\}_{a,b}$ form an open basis of the constructible topology and the~$\{C(a,b)\compl\}_{a,b}$ form a closed basis.
A \emph{proconstructible} (or \emph{patch-closed}) subset is a closed subset for the constructible topology, \ie an arbitrary intersection of subsets~$C(a,b)\compl=\supp(a)\cup \supp(b)\compl=\SET{\cP\in\SpcK}{a\notin\cP\textrm{ or }b\in\cP}$.
\end{rec}

\begin{rec}
\label{rec:Thom-and-proconst}%
A subset~$Y\subseteq \SpcK$ is proconstructible and Thomason (if and) only if $Y$ is closed with~$Y\compl$ quasi-compact (see \cite[Theorem 1.5.4\,(i) and Corollary 4.1.6\,(ii)]{DickmannSchwartzTressl19}).
\end{rec}

\begin{prop}
\label{prop:spc(AI)-procon}%
Let $\cI\subseteq \mK$ be a Serre $\otimes$-ideal. Then $\WI$ is a proconstructible subset of~$\SpcK$; more precisely, if we let  $\cM(\cI)\coloneq\SET{(a,b)\in\cK\times\cK}{\bar{b}\in\ideal{\bar{a}}\textrm{ in }\AI}$ then
\begin{equation}
\label{eq:spc(AI)-procon}%
\WI=\bigcap_{(a,b)\in \cM(\cI)}\supp(a)\cup \supp(b)\compl.
\end{equation}
\end{prop}
\begin{proof}
Let $\cH\supseteq\cI$ be a homological prime in~$\WhI$ and let $(a,b)\in \cM(\cI)$, so that $\bar{b}$ belongs to the Serre $\otimes$-ideal~$\ideal{\bar{a}}$ generated by~$\bar{a}$ in~$\AI$; a fortiori $\bar{b}\in\ideal{\bar{a}}$ in the further quotient~$\AH$. This relation prevents $\bar{a}=0$ and $\bar{b}\neq0$ in~$\AH$, that is, we cannot have $a\in\yon\inv(\cH)=\nos(\cH)$ and $b\notin\yon\inv(\cH)=\nos(\cH)$; hence we have $\nos(\cH)\in \supp(a)\cup \supp(b)\compl$.
This proves~$\subseteq$ in~\eqref{eq:spc(AI)-procon}.

Conversely, let $\cP\in\SpcK$ belong to the right-hand side of~\eqref{eq:spc(AI)-procon}.
Define $\cJ=\ideal{\cI\,\cup\,\yon(\cP)}$ the Serre $\otimes$-ideal of~$\mK$ generated by~$\cI$ and~$\cP$. Let $S$ be the $\otimes$-multiplicative class~$\cK\sminus\cP$.
We claim that $\cJ\cap \yon(S)=\varnothing$. \textsl{Ab absurdo}, if $b\in\cK\sminus \cP$ satisfies $\yon(b)\in\cJ$ then there exist~$a\in\cP$ such that $\yon(b)\in \ideal{\cI\cup\{\yon(a)\}}$. Modulo~$\cI$ we get $\bar{b}\in \ideal{\bar{a}}$ in~$\AI$, that is, $(a,b)\in \cM(\cI)$. By assumption $\cP\in \supp(a)\cup \supp(b)\compl$ which gives us $a\notin\cP$ or~$b\in \cP$, either of which is a contradiction. So $\cJ\cap \yon(S)=\varnothing$ and by~\cite[Lemma~3.8]{Balmer20a}, there exists~$\cH\in\SpchK$ such that $\cJ\subseteq \cH$ and~$\cH\cap \yon(S)=\varnothing$. By construction $\cI\subseteq \cJ\subseteq \cH$, so $\cH\in \WhI$ and its image $\nos(\cH)=\yon\inv(\cH)$ in~$\SpcK$ contains~$\cP$ and avoids~$S=\cK\sminus\cP$, hence equals~$\cP$.
\end{proof}

\begin{exas}
Let us discuss subsets of $\SpcK$ that cannot be the spectrum of a layer, because they are not proconstructible. We use~\Cref{rec:Thom-and-proconst}.
Suppose that $\SpcK$ has a generic point~$\eta$ such that~$\{\eta\}$ is not open; for instance $\cK$ could be~$\SH\compact$ (see \cite[Corollary 9.5(c)]{Balmer10b}), or~$\Dperf(R)$ for $R$ a domain. Then $\SpcK \sminus \{\eta\}$ is not the spectrum of any layer.
\end{exas}

\begin{rem}
We want to prove the converse to \Cref{prop:spc(AI)-procon}: Every proconstructible is the spectrum~$\WI$ of a layer. In general, there are many such~$\cI$ and we describe the \emph{largest} one in \Cref{prop:Imax}. But we want to use the same quantifiers as for a proconstructible subset: If we {choose} $\cN\subseteq \cK\times \cK$ then $\bigcap_{(a,b)\in\cN}\supp(a)\cup\supp(b)\compl$ is proconstructible by \Cref{rec:constructible}. We want an~$\cI$ for every such choice of~$\cN$.
We need a preparation.
\end{rem}

\begin{lem}
\label{lem:h(a)-or-ker(h(coev))-in-H}
For any $\cH \in \Spch(\cK)$ and $a \in \cK$, exactly one of the following is true:
\begin{enumerate}[\rm(i)]
\item
\label{it:h(a)-H}
The object $\yon(a)$ belongs to~$\cH$.
\item
\label{it:ker(h(coev))-H}%
The kernel $\ker(\yon(\coev_a\colon \unit\to a^\vee\otimes a))$ of coevaluation belongs to~$\cH$.
\end{enumerate}
\end{lem}

\begin{proof}
If~$\bar{a}=0$ in~$\AH$ then~\eqref{it:h(a)-H}. If~$\bar{a}\neq 0$ then by \Cref{rec:dual} the kernel of $\coev_{\bar{a}}\colon \bar{\unit}\to \bar{a}^\vee\otimes\bar{a}$ in~$\AH\fp$ belongs to the Serre $\otimes$-ideal~$\Ker(\bar{a}\otimes-)\subsetneqq\AH\fp$ which is zero; this gives~\eqref{it:ker(h(coev))-H}.
\end{proof}

\begin{lem}
\label{lem:claim}%
Let $a,b\in\cK$ and~$\cH\in\SpchK$. Then $\ker(\yon(\coev_a)) \otimes \yon(b)$ belongs to~$\cH$ if and only if~$\nos(\cH)\in\supp(a)\cup\supp(b)\compl$.
\end{lem}
\begin{proof}
Since $\yon(b)$ is~$\otimes$-flat and~$\cH$ is a maximal Serre~$\otimes$-ideal we have $\ker(\yon(\coev_a)) \otimes \yon(b)\in\cH$ if and only if~$\ker(\yon(\coev_a))\in\cH$ or~$\yon(b)\in\cH$. By~\Cref{lem:h(a)-or-ker(h(coev))-in-H} this is equivalent to~$\yon(a)\notin\cH$ or~$\yon(b)\in\cH$, which reads~$\nos(\cH)\in\supp(a)\cup\supp(b)\compl$.
\end{proof}

\begin{prop}
\label{prop:procon-Spc(AI)}%
Every proconstructible subset $C$ of~$\SpcK$ is of the form~$\WI$ for a Serre $\otimes$-ideal~$\cI$ of~$\mK$: For any choice of~$\cN\subseteq \cK\times \cK$ such that $C=\bigcap_{(a,b)\in \cN}\supp(a)\cup\supp(b)\compl$, define a Serre $\otimes$-ideal of~$\mK$ by the following generators
\begin{equation*}
\cI =\cI(\cN)\coloneq \ideal{ \SET{ \ker(\yon(\coev_a))\otimes \yon(b)  }{ (a,b) \in \cN}}.
\end{equation*}
Then we have $\WhI=\nos\inv(C)$ and~$\WI=C$.
\end{prop}
\begin{proof}
Let~$\cH\in\SpchK$. We have $\cH\in\WhI$, \ie $\cI \subseteq \cH$ if and only if $\ker(\yon(\coev_a)) \otimes \yon(b) \in \cH$ for every $(a, b) \in \cN$. By~\Cref{lem:claim}, this is equivalent to $\nos(\cH)\in\supp(a)\cup \supp(b)\compl$ for every~$(a,b)\in\cN$, or in other words~$\nos(\cH)\in \bigcap_{(a,b)\in \cN}\supp(a)\cup\supp(b)\compl=C$.
So we have proved~$\WhI=\nos\inv(C)$, which implies~$\WI=\nos(\nos\inv(C))=C$ by surjectivity of~$\nos$.
\end{proof}
\begin{rem}
In \Cref{prop:procon-Spc(AI)}, one can take for~$\cN$ the largest one, namely $\cN=\cN(C) = \SET{ (a, b) \in \cK \times \cK }{ C \subseteq \supp(a)\cup \supp(b)\compl}$. See \Cref{rec:constructible}.
\end{rem}

Let us discuss the special case of closed subsets (that are in particular proconstructible).

\begin{dfn}
\label{dfn:I_S}%
Let $S\subseteq \cK$ be a $\otimes$-multiplicative class of objects in $\cK$. The set $Z(S) \coloneq \bigcap_{a \in S} \supp(a)$ is closed in $\Spc(\cK)$; in fact any closed subset of $\SpcK$ is of this form. Using the notation of \Cref{prop:procon-Spc(AI)}, we set $\cN_S = \{(a, \unit) \mid a \in S \} \subseteq \cK \times \cK$ and we define
    \begin{equation*}
        \cI_S \coloneq \cI(\cN_S) = \ideal{ \{  \ker(\yon(\coev_a)) \mid a \in S\}}.
    \end{equation*}
\end{dfn}
\begin{cor}
\label{cor:I_S}%
We have $\W(\cI_S) = Z(S)=\cap_{a\in S}\supp(a)$ and~$\Wh(\cI_S)=\nos\inv(Z(S))$.
\end{cor}
\begin{proof}
Direct from \Cref{prop:procon-Spc(AI)} for~$\cN=\cN_S$ since $\supp(\unit)\compl=\varnothing$.
\end{proof}

\begin{prop}
\label{prop:I_S}%
Let $S$ be a $\otimes$-multiplicative subset of objects in $\cK$. Then we have 
\begin{equation}
\label{eq:I_S-equality}
\cI_{S} = \SET{M \in \mK}{\textrm{there exists } a \in S\textrm{ such that }\yon(a) \otimes M = 0}.
\end{equation}
\end{prop}
\begin{proof}
Since $S$ is $\otimes$-multiplicative, the right-hand side of \eqref{eq:I_S-equality} is a Serre $\otimes$-ideal of~$\mK$. Take $a\in S$. By \Cref{rec:dual} we have $\yon(a) \otimes \ker(\yon(\coev_a))=0$. This proves the inclusion~$\subseteq$ in~\eqref{eq:I_S-equality}.
Conversely, let~$M \in \mK$ and $a \in S$ with~$\yon(a) \otimes M =0$. Let~$\iota\colon M \into \yon(c)$ be a monomorphism for some~$c \in \cK$ and consider the commutative square in~$\mK$
\[\begin{tikzcd}[row sep=1.5em]
	M && {\yon(c)} \\
	{M \otimes \yon(a \otimes a^{\vee})} && {\yon(c) \otimes \yon(a \otimes a^{\vee})}
	\arrow["\iota", tail, from=1-1, to=1-3]
	\arrow["{M \otimes \yon(\coev_a)}"', from=1-1, to=2-1]
	\arrow["{\yon(c) \otimes \yon(\coev_a)}", from=1-3, to=2-3]
	\arrow["{\iota \otimes \yon(a \otimes a^{\vee})}"', from=2-1, to=2-3]
\end{tikzcd}
\]
In the quotient~$(\mK)/\ideal{\ker(\yon(\coev_a))} $, the morphism $\yon(c) \otimes \yon(\coev_a)$ becomes a monomorphism, hence so does $M \otimes \yon(\coev_a)$. Thus $M \otimes \yon(a \otimes a^{\vee}) = 0$ forces $M \in \ideal{\ker(\yon(\coev_a))} $.
\end{proof}

\begin{rem}
\label{rem:I_S-Thomason}%
With the notation of \Cref{prop:I_S}, suppose that the closed subset $Z(S)=\cap_{a\in S}\supp(a)$ is Thomason. Then $Z\compl=\cup_{a\in S}\supp(a)\compl$ is quasi-compact, hence there exists $a_0\in S$ such that $Z(S)=\supp(a_0)$.
In that case we have the following modification of~\eqref{eq:I_S-equality}:
\begin{equation*}
\label{eq:I_S-newequality}
\cI_{S} = \SET{M \in \mK}{\textrm{there exists } a \in S\textrm{ with } \supp(a) = Z(S) \textrm{ such that }\yon(a) \otimes M = 0}.
\end{equation*}
Indeed, if $\yon(b)\otimes M=0$ for some~$b\in S$ then $\yon(a_0\otimes b)\otimes M=0$ and $\supp(a_0\otimes b)=Z(S)$.
\end{rem}

To wrap up the section, we describe the largest Serre $\otimes$-ideal with a given support. It is built in the same way as the radical of an ideal in commutative algebra.

\begin{dfn}
\label{dfn:Imax}%
For~$C\subseteq\SpcK$ we set~$\Imax(C)=\bigcap_{\cH\in\nos\inv(C)}\cH$. It is a Serre $\otimes$-ideal.
\end{dfn}
\begin{prop}
\label{prop:Imax}%
Let~$C\subseteq\SpcK$ be a subset and~$\cI\subseteq \mK$ a Serre $\otimes$-ideal.
We have $\WhI\supseteq \nos\inv(C)$ if and only if~$\cI\subseteq\Imax(C)$.
\end{prop}

\begin{proof}
If $\WhI\supseteq \nos\inv(C)$ then $\cI \subseteq \bigcap_{\cH \in \WhI} \cH \subseteq \bigcap_{\cH \in \nos\inv(C)} \cH =\Imax(C)$. Conversely, if $\cI\subseteq\Imax(C)$ then for all $\cH \in \nos\inv(C)$ we have~$\cH \supseteq\Imax(C) \supseteq \cI$. Therefore $\cH \in \WhI$. 
\end{proof}

As announced, $\Imax(C)$ is another solution to \Cref{prop:procon-Spc(AI)}.
\begin{prop}
\label{prop:Spc(Imax)}%
Let $C\subseteq\SpcK$ be a proconstructible subset. Then $\Wh(\Imax(C))=\nos\inv(C)$ and consequently $\W(\Imax(C))=C$.
\end{prop}
\begin{proof}
Write~$C=\bigcap_{(a,b)\in\cN}\supp(a)\cup\supp(b)\compl$ for some subset~$\cN\subseteq\cK\times\cK$.
The inclusion $\nos\inv(C) \subseteq \Wh(\Imax(C))$ holds by~\Cref{prop:Imax} for~$\cI=\Imax(C)$.
Conversely, recall~$\cI(\cN)=\ideal{\SET{\ker(\coev_{\yon(a)})\otimes\yon(b)}{(a,b)\in\cN}}$ from \Cref{prop:procon-Spc(AI)}. 
Since the generators of~$\cI(N)$ belong to every~$\cH\in \nos\inv(C)$ by \Cref{lem:claim}, we have $\cI(\cN)\subseteq\Imax(C)$.
We deduce from this that $\Wh(\Imax(C))\subseteq \Wh(\cI(N))=\nos\inv(C)$ using \Cref{prop:procon-Spc(AI)}.
\end{proof}

\section{Layers over Thomason subsets}
\label{sec:Thom-layers}%

We keep~$\cK$ as in \Cref{hyp:K-T}.
In \Cref{sec:procon} we have seen that the spectrum~$\WI$ of every~layer~$\AI$ is a proconstructible subset of~$\SpcK$. We now want to study the layers that live over a Thomason subset~$Y\subseteq \SpcK$ by using the tt-ideal~$\cK_Y=\SET{c\in\cK}{\supp(c)\subseteq Y}$.

\begin{lem}
\label{lem:ker(i*)}%
Let $Y \subseteq \SpcK$ be a Thomason subset.
Let $i^{\ast}\colon \MK \to \MK_{Y}$ denote the restriction along the inclusion $i \colon \cK_{Y} \hook \cK$. Then $i^{\ast}$ is a Gabriel quotient with kernel
\begin{equation}
\label{eq:ker(i*)}%
    \Ker(i^{\ast}) = \SET{ M \in \MK}{\yon(a) \otimes M = 0 \textrm{ for all }a\in \cK_Y}.
\end{equation}
Thus $\Ker(i^{\ast})$ is a $\otimes$-ideal, and $\yon(a) \otimes N$ is $\Ker(i^{\ast})$-local for all~$a \in \cK_Y$ and~$N \in \MK$.  
\end{lem}
\begin{proof}
Since $i \colon \cK_{Y} \hook \cK$ is fully faithful, the induced $i_{!}\colon \MK_Y\to \MK$ is fully faithful. Hence its exact right adjoint~$i^{\ast}\colon \MK\to \MK_Y$ is a Gabriel quotient.
To prove~\eqref{eq:ker(i*)}, recall that, by \cite[Lemma 2.5]{BalmerKrauseStevenson19}, there is a natural isomorphism of functors $\cK^{\rm{op}} \to \rm{Ab}$
\begin{equation}\label{eq:tensor-with-representables}
    \yon(a) \otimes M \cong M(a^{\vee} \otimes -)
\end{equation}
for any~$M \in \MK$ and~$a \in \cK$. 
Hence if~$i^{\ast}(M)=M\restr{\cK_Y}=0$ then for every~$a\in\cK_Y$ we have $\yon(a)\otimes M=0$ since $a^\vee\otimes\cK\subseteq\cK_Y$. Conversely, if $\yon(a)\otimes M=M(a^\vee\otimes -)=0$ then $M$ vanishes on $a\otimes a^\vee\otimes a$ and therefore on~$a$ itself, which is a direct summand.
Hence~\eqref{eq:ker(i*)}. It follows that $\Ker(i^*)$ is a $\otimes$-ideal of~$\MK$. Finally, fix $a \in \cK_Y$ and $N \in \MK$. Then $(- \otimes \yon(a^{\vee})) \,\adj\, (\yon(a) \otimes -)$ is an adjunction of exact functors on~$\cA=\MK$. Therefore
\begin{equation*}
    \Ext^*_{\cA}\big(M, \yon(a) \otimes N\big) \cong \Ext^*_{\cA}\big( M \otimes \yon(a^{\vee}), N\big) = 0
\end{equation*}
for every~$M \in \Ker(i^{\ast})$.
This implies that $\yon(a) \otimes N$ is $\Ker(i^{\ast})$-local.
\end{proof}
The localizing subcategory~$\Ker(i^*)$ is almost never generated by finitely presented objects:
\begin{prop}
\label{prop:Bmin-not-finite-type}%
With the same notation as in \Cref{lem:ker(i*)}, the localizing subcategory~$\Ker(i^*)$ is of finite type if and only if $Y$ is closed and open, in which case~$\cK \cong \cK_Y \times \cK_{Y\compl}$.
\end{prop}
\begin{proof}
By \Cref{thm:finite-type}, $\Ker(i^*)$ is of finite type if and only if $i^*(\unit)$ is finitely presented in~$\MK_Y$. So suppose first that $i^*(\unit) \in \mK_Y$. Note that $i^*(\unit)$ remains a cohomological functor, so it is injective in $\mK_Y$ by \cite[Lemma 2.7]{Krause00}. Hence $i^*(\unit)\simeq \yon(e)$ for some $e \in \cK_Y$ (\Cref{rec:Yoneda}). It follows (see \Cref{rec:eY-fY}) that $\eY=i_! i^*(\unit) \simeq i_!(\yon(e)) \cong \yon(i(e))$ is finitely presented in~$\cA$, and therefore so is~$\fY$ by~\eqref{eq:LES-eY-fY}. Since $\fY$ is also cohomological, we get $\fY \cong \yon(f)$ for some $f \in \cK$ by~\cite[Lemma 2.7]{Krause00} again.
Then $\Hom(\Sigma\inv \fY, \eY) = \Hom(\fY, \Sigma \eY) = 0$ by $\eY\otimes\fY=0$ \eqref{eq:fY-tens-KY=0} and rigidity. Thus there is a short exact sequence $0 \to \yon(e) \to \unit \to \yon(f) \to 0.$ It easily follows that $\unit \simeq e \oplus f$ in~$\cK$, forcing $\cK \cong \cK_Y \times \cK_{Y\compl}$ with~$Y\compl=\supp(f)$ closed.
The reverse direction is clear: $i^*(\unit)\cong \yon(\unit_{\cK_Y})$.
\end{proof}

\begin{dfn}
\label{dfn:Bmin}%
Let~$Y\subseteq\SpcK$ be Thomason. We denote by~$\Bmin(Y)$ the localizing $\otimes$-ideal of $\MK$ of \Cref{lem:ker(i*)}
(the notation~$\Bmin$ is justified by \Cref{thm:Bmin} and \Cref{cor:bookend}):
\begin{equation*}
\Bmin(Y)\coloneq\Ker(\MK\overset{i^*}{\,\onto\,} \MK_Y)
= \SET{ M \in \MK}{\yon(a) \otimes M = 0 \textrm{ for all }a\in \cK_Y}.
\end{equation*}
%
\end{dfn}

\begin{thm}
\label{thm:Bmin}%
Let~$Y\subseteq\SpcK$ be a Thomason subset and~$\cI\subseteq \mK$ a Serre $\otimes$-ideal.
We have $\WhI\subseteq \nos\inv(Y)$ if and only if~$\Bmin(Y)\subseteq \Loc(\cI)$.
\end{thm}

\begin{proof}
First, suppose that $\WhI\subseteq \nos\inv(Y)$. For every homological prime~$\cH \in \nos \inv(Y)$, we have $\overline{\fY} =0$ in~$\AH$ by~\Cref{lem:eY-fY-H}\,\eqref{it:eY-fY-H-in}. So $\fY$ is zero in the residue fields $\AH$ for all~$\cH\in\WhI$.
As $\fY$ is a $\otimes$-flat ring object, it must go to zero in~$\AI$ by the Nilpotence Theorem (\Cref{rec:nilpotence}).
Therefore $\overline{\eY}\simeq \bar\unit$ in~$\AI$ by~\eqref{eq:LES-eY-fY}. Let $M\in \Bmin(Y)$. By~\eqref{eq:ker(i*)} we have $M \otimes \yon(a) = 0$ for every $a\in \cK_Y$. Since $\eY$ is a colimit of~$\yon(a)$ for $ a \in \cK_Y$, we have $M\otimes \eY=0$ as well. Therefore in~$\AI$, we have $\bar M\cong \bar M\otimes \bar \unit \simeq \bar M\otimes \overline{\eY}= 0$, proving that $M\in \Loc(\cI)$. 

Conversely, suppose that $\WhI\not\subseteq \nos\inv(Y)$. Then there exists~$\cH\supseteq\cI$ with $\nos(\cH)\notin Y$. We have $\kappa(\cH)\notin\Loc(\cI)$ by~\Cref{lem:H-in-Spch(A/I)}. On the other hand, because $\nos(\cH) \notin Y$, we must have $\yon(a) \in \cH$ for all $a \in \cK_Y$.
By the projection formula for $Q_{\cH} \adj R_{\cH}$, for all $a \in \cK_Y$, we get in~$\cA$
\begin{equation*}
    \kappa(\cH)\otimes\yon(a) = R_{\cH}\left(\overline{\kappa(\cH)}\right) \otimes \yon(a) \cong R_{\cH} \left( \overline{\kappa(\cH)} \otimes  Q_{\cH}(\yon(a) ) \right) =0.
\end{equation*}
Therefore by~\eqref{eq:ker(i*)}, we have~$\kappa(\cH) \in \Bmin(Y)$. Hence we have proved that $\Bmin(Y)$ is not contained in~$\Loc(\cI)$ since the former contains~$\kappa(\cH)$ but the latter does not.
\end{proof}

We can combine the above with \Cref{prop:Imax} when $Y$ is furthermore proconstructible (see~\Cref{rec:Thom-and-proconst}). We obtain an `interval' of Serre $\otimes$-ideals whose spectrum is~$Y$.
\begin{cor}
\label{cor:bookend}%
Let~$Y\subseteq\SpcK$ be a closed subset with quasi-compact complement and let~$\cI\subseteq \mK$ be a Serre $\otimes$-ideal. Consider the two localizing $\otimes$-ideals of~$\cA=\MK$
\[
\Bmin(Y)=\Ker(i^*\colon \MK\onto \MK_Y)
\qquadtext{and}
\Bmax(Y)\coloneq\Loc(\Imax(Y)).
\]
Then we have $\WhI=\nos\inv(Y)$ if and only if~$\Bmin(Y)\subseteq \Loc(\cI)\subseteq\Bmax(Y)$.
\end{cor}

\begin{proof}
By~\Cref{prop:Imax},~$\WhI\supseteq\nos\inv(Y)$ if and only if
$\cI\subseteq \Imax(Y)$, that is, $\Loc(\cI)\subseteq \Bmax(Y)$.
By \Cref{thm:Bmin}, $\WhI\subseteq\nos\inv(Y)$ if and only if~$\Bmin(Y)\subseteq\Loc(\cI)$.
\end{proof}

\begin{rem}
Although $\Bmax(Y)$ is of finite type, the subcategory $\Bmin(Y)$ need not be (\Cref{prop:Bmin-not-finite-type}). This may explain why $\Bmin(Y)$ does not have an analog in commutative algebra.
\end{rem}

\section{Preservation of layers}
\label{sec:Y-invariance}%

For the whole section,~$F\colon \cK\to \cL$ is a tt-functor between two rigid idempotent-complete tt-categories $\cK$ and~$\cL$, as in~\Cref{hyp:K-T}. We denote the induced map on spectra by
\begin{equation*}
\varphi=\Spc(F)\colon \SpcL\to \SpcK.
\end{equation*}
Recall that $\supp_{\cL}(F(c))=\varphi\inv(\supp_{\cK}(c))$ for all~$c\in\cK$.
It follows that the preimage $\varphi\inv(Y)$ of a Thomason subset~$Y$ is Thomason and that $F(\cK_Y)\subseteq \cL_{\varphi\inv(Y)}$.

\begin{dfn}
\label{dfn:Y-equivalence}%
Let $Y\subseteq\SpcK$ be a Thomason subset. We say that~$F$ is an \emph{equivalence over~$Y$} if the restriction of $F$ on~$\cK_Y$ is an equivalence
\begin{equation*}
F_Y\coloneq F\restr{\cK_Y}\colon \cK_Y\isotoo \cL_{\varphi\inv(Y)}.
\end{equation*}
(One might prefer to call such an~$F$ an `equivalence \emph{with supports in}~$Y$' since $\cK_Y$ should not be mistaken for an honest tt-category.
In particular, by~\Cref{rem:not-T_Z}, if~$Y\hook X$ is a closed subscheme and~$\cK=\Dperf(X)$ then $\cK_Y$ should not be confused with~$\Dperf(Y)$.)
\end{dfn}

The following result also holds for~$\Spch$; see \Cref{cor:Y-equiv-Spch}.
\begin{prop}
\label{prop:Y-equivalence}%
Let $F\colon \cK\to \cL$ be a tt-functor that is an equivalence over~$Y$. Then $\varphi$ restricts to a continuous bijection~$\varphi\inv(Y)\isoto Y$. Moreover, if~$Z$ is a Thomason subset of~$\SpcL$ contained in~$\varphi\inv(Y)$ then $\varphi(Z)$ is a Thomason subset of~$\SpcK$. Finally, if~$Y$ is closed with quasi-compact complement then this bijection~$\varphi\inv(Y)\isoto Y$ is a homeomorphism.
\end{prop}
\begin{proof}
Write $Y'=\varphi\inv(Y)$ and consider $\varphi\restr{Y'}\colon Y'\to Y$. Since $F$ detects nilpotence on~$\cK_Y$ it follows from~\cite{Balmer18a} that~$\varphi\restr{Y'}$ is surjective. 
For injectivity, we adapt~\cite[Lemma~4.5]{BalmerSanders25}. Let $\cQ,\cQ'\in Y'$ satisfy~$\varphi(\cQ)=\varphi(\cQ')\in Y$. Let $c\in\cK_Y$ such that $\varphi(\cQ)=\varphi(\cQ')\in \supp(c)$. This implies that $F(c)\notin \cQ^{'}$. Let $d\in\cQ$. Then we have $d\otimes F(c)\in\cL_{Y'}$. By our equivalence~$\cK_Y\isoto \cL_{Y'}$, there exists $b\in\cK_Y$ such that $F(b)\simeq d\otimes F(c)\in\cQ$. Hence $b\in F\inv(\cQ)=F\inv(\cQ')$ and therefore $d\otimes F(c)\simeq F(b)\in\cQ'$. Since $F(c)\notin\cQ'$ we have $d\in\cQ'$. In short $\cQ\subseteq \cQ'$. By symmetry, $\cQ=\cQ'$. So $\varphi\restr{Y'}$ is a continuous bijection.

Let $d\in \cL_{Y'}$ and~$c\in\cK_{Y}$ such that $F(c)\simeq d$. Then~$\varphi\inv(\supp(c))=\supp(d)$ implies that $\supp(c)=\varphi(\supp(d))$. A Thomason subset~$Z\subseteq\SpcL$ contained in~$Y'$ is a union of~$\supp(d)$ with~$d\in\cL_{Y'}$. This gives the `moreover part'. When $Y$, hence~$Y'$, is closed Thomason (\ie the support of an object), then these~$\{\supp(d)\}_{d\in\cL_{Y'}}$ form a basis of the subspace topology on~$Y'$. Hence the final statement.
\end{proof}

We want to show that equivalences over Thomason subsets localize well.
\begin{rec}
\label{rec:K(U)}%
If~$U\subseteq \SpcK$ is the complement of a Thomason subset (\ie an intersection of quasi-compact open subsets), then~$\cK(U)$ denotes~$(\cK/\cK_{U\compl})^\natural$ with the obvious functor $\res_U\colon \cK\onto \cK(U)$. One has $\Spc(\cK(U))\cong U$ and~$\supp(\res_U(c))=\supp(c)\cap U$ for all~$c\in\cK$.
\end{rec}

\begin{prop}
\label{prop:excision}%
Let $F\colon \cK\to \cL$ be a tt-functor that is an equivalence over a Thomason~$Y$. 
Let $U\subseteq \SpcK$ and $V\subseteq\SpcL$ with~$U\compl$ and~$V\compl$ Thomason. Suppose that $V\subseteq\varphi\inv(U)$ so that~$F$ restricts to~$\bar F\colon \cK(U)\to \cL(V)$. If~$\bar F$ is conservative then it is an equivalence over~$Y\cap U$.
\end{prop}
\begin{proof}
Since we assume $\varphi\inv(U\compl)\subseteq V\compl$ we have $F(\cK_{U\compl})\subseteq\cL_{V\compl}$ and~$F$ descends to a tt-functor~$\bar F\colon \cK(U)\to \cL(V)$ such that $\res_V \circ F\cong\bar F\circ \res_U$.
Let $Y'=\varphi\inv(Y)$. We have
\begin{equation}
\label{eq:excision}%
\cK(U)_{Y\cap U}\cong (\cK_Y/\cK_{Y\cap U\compl})^\natural
\qquadtext{and}
\cL(V)_{Y'\cap V}\cong(\cL_{Y'}/\cL_{Y'\cap V\compl})^\natural.
\end{equation}
Indeed, if $Z:=U\compl$ then $\cK_{Y\cap Z}=\cK_Y\cap \cK_Z$ and the functor~$(\cK_{Y}/\cK_{Y\cap Z})^\natural\to (\cK/\cK_Z)^\natural=\cK(U)$ induced by~$\res_U$ is fully faithful by~\cite[Chapitre II, Proposition 2.3.5]{Verdier96}. Its image is a tt-ideal determined by its support, which is obviously contained in~$Y\cap U$ since $\supp(\res_U(c))=\supp(c)\cap U$. Conversely, for every $d\in \cK(U)$ there exists~$c\in \cK$ such that $\res_U(c)\simeq d\oplus \Sigma(d)$; if $d\in\cK(U)_{Y\cap U}$ we get $\supp(c)\subseteq\supp(d)\cup Z\subseteq Y\cup Z$. Since~$\cK_{Y\cup Z}=\cK_Y\vee\cK_Z$ has the same image in~$\cK(U)$ as~$\cK_Y$, we get the first equivalence in~\eqref{eq:excision}. The second is the same in~$\cL$. We claim that we have a commutative diagram as follows (where ${}\restr{\cdots}$ indicates restriction)
\[
\xymatrix@C=4em@R=1.5em{
\cK \ar[d]_-{F} \ar@{->>}[r]_-{\res_U}
& \cK(U) \ar[d]_-{\bar F}
& \cK(U)_{Y\cap U} \ar@{ >->}[l]^-{\incl} \ar[d]_-{\bar F\restr{\cdots}}
& \cK_Y \ar@{->>}[l]^-{\res_U} \ar[d]_-{F\restr{\cdots}}^-{\simeq} \ar@{ >->}@/_2em/[lll]^-{\incl}
& \cK_{Y\cap U\compl} \ar[d]_-{F\restr{\cdots}}^-{\simeq} \ar@{ >->}[l]^-{\incl}
\\
\cL \ar@{->>}[r]^-{\res_V}
& \cL(V)
& \cL(V)_{Y'\cap V} \ar@{ >->}[l]_-{\incl}
& \cL_{Y'} \ar@{->>}[l]_-{\res_V} \ar@{ >->}@/^2em/[lll]_-{\incl}
& \cL_{Y'\cap V\compl} \ar@{ >->}[l]_-{\incl}
}
\]
The functor $\bar{F}$ sends~$\cK(U)_{Y\cap U}$ into~$\cL(V)_{Y'\cap V}$ thanks to~\eqref{eq:excision} and the fact that for every~$c\in\cK_Y$, we have $\supp_{\cL(V)}(\bar F(\res_U(c)))=\supp_{\cL(V)}(\res_V(F(c)))=\varphi\inv(\supp(c))\cap V\subseteq Y'\cap V$.
In the right-most square, the equivalence $F_Y:=F\restr{\cK_Y}\colon\cK_Y\isoto \cL_{Y'}$ sends the subcategory~$\cK_{Y\cap U\compl}$ exactly onto~$\cL_{Y'\cap V\compl}$ for these subcategories are respectively $\Ker(\res_U)\cap\cK_Y$ and $\Ker(\res_V)\cap\cL_{Y'}$ by~\eqref{eq:excision} again, and $0=\res_V(F_Y(c))=\bar{F}(\res_U(c))$ if and only if~$\res_U(c)=0$ by our conservativity hypothesis on~$\bar{F}$. It follows from~\eqref{eq:excision} that $F_Y$ descends to an equivalence~$\cK(U)_{Y\cap U}\isoto \cL(V)_{Y'\cap V}$, which is necessarily~$\bar{F}\restr{\cK(U)_{Y\cap U}}$ by the above diagram.
\end{proof}

\begin{lem}
\label{lem:equiv-on-Y-Hom}%
Let $F\colon \cK\to \cL$ be a tt-functor that is an equivalence over~$Y$. 
Then~$F$ yields a natural isomorphism~$\alpha_{b,c}\colon \Homcat{K}(b,c)\isoto \Homcat{L}(F(b),F(c))$ for every $b\in \cK_Y$ and~$c\in\cK$.
\end{lem}
\begin{proof}
Naturality is easy.
Then, since every $b\in\cK$ is a retract of~$b\otimes b^\vee\otimes b$, it suffices to prove that~$\alpha_{b\otimes b^\vee\otimes b\,,\,c}$ is an isomorphism.
One verifies that~$\alpha_{b\otimes b^\vee \otimes b,c}$ factors as ${\Homcat{K}(b\otimes b^\vee\otimes b, c)}\cong{\Homcat{K}(b \otimes b,b\otimes c)}
\underset{\alpha_{b\otimes b,b\otimes c}}{\isoto}
{\Homcat{L}(F(b \otimes b),F(b\otimes c))} \cong {\Homcat{L}(F(b\otimes b^\vee\otimes b), F(c))}$
where the unlabelled isomorphisms are given by rigidity and $\alpha_{b\otimes b,b\otimes c}$ is an isomorphism since~$b\otimes c\in\cK_Y$ and $F\restr{\cK_Y}$ is fully faithful by assumption. Hence the result.
\end{proof}
\begin{lem}
\label{lem:F!-restr}%
Let $F\colon \cK\to \cL$ be a tt-functor that is an equivalence over~$Y$. Let $Y' = \varphi\inv(Y) \subseteq \SpcL$ be the preimage of~$Y$.
Then the following diagram commutes up to isomorphism:
\begin{equation}
\label{eq:F!-restriction}%
\vcenter{\xymatrix@C=3em@R=1.5em{
\MK \ar[r]^-{F_!} \ar[d]_-{i^*}
& \ML \ar[d]^-{j^*}
\\
{\MK_Y} \ar[r]^-{(F_Y)_!}_-{\simeq}
& {\ML_{Y'}}
}}
\end{equation}
where $i\colon \cK_Y\hook\cK$ and~$j\colon \cL_{Y'}\hook\cL$ are the inclusions and where $F_Y=F\restr{\cK_Y}\colon \cK_Y\isoto \cL_{Y'}$.
\end{lem}
\begin{proof}
This is not trivial since we mix~$(-)^*$ and~$(-)_!$. As $F_Y\colon \cK_Y\isoto \cL_{Y'}$ is an equivalence, so is~$(F_Y)_!$, and its inverse is its adjoint~$(F_Y)^*$. Hence it suffices to construct an isomorphism
\[
i^*\isoto (F_Y)^*\circ j^*\circ F_! \colon \MK\to \MK_Y.
\]
All four functors are cocontinuous and exact. By the universal property of~$\MK$ (\cite[Lemma~2.2]{Krause00}), it suffices to give such an isomorphism after precomposing with the Yoneda embedding $\yon\colon \cK\hook \MK$. For every~$c\in\cK$, we have a natural transformation $\alpha$ in~$\MK_Y$
\[
\alpha_{-,c}\colon i^*(\yon(c))=\Homcat{K}(-,c)\restr{\cK_Y}\overset{F}{\,\longrightarrow\,}
\Homcat{\cL}(F(-),F(c))\restr{\cK_Y}=(F_Y)^*(j^*(F_!(\yon(c)))).
\]
So it suffices to know that for all~$c\in\cK$ and all~$b\in \cK_Y$ the functor~$F$ yields a bijection $\alpha_{b,c}\colon\Homcat{K}(b,c)\isoto \Homcat{L}(F(b),F(c))$, which is what we proved in~\Cref{lem:equiv-on-Y-Hom}.
\end{proof}

Recall $\WhI=\SET{\cH\in\SpchK}{\cI\subseteq \cH}$ and~$\WI=\nos(\WhI)$ from~\Cref{dfn:spec-layer} and recall that $\nos\colon \Spch\onto \Spc$ is often a bijection in examples (\Cref{rec:Spch}). Therefore, in first approximation, the reader can ignore~$\nos$ below and think~$\WhI=\WI$.
\begin{thm}
\label{thm:main}%
Let $F\colon \cK\to \cL$ be a tt-functor and~$\varphi=\Spc(F)$. Let $Y\subseteq \SpcK$ be a Thomason subset such that $F$ is an equivalence over~$Y$ (\Cref{dfn:Y-equivalence}). Let $Y'=\varphi\inv(Y)$. Then, using $F_!\colon \mK\to \mL$, we get an inclusion-preserving bijection
\[
\left\{\ 
\substack{\displaystyle\textrm{Serre $\otimes$-ideals~$\cI$ of~$\mK$}
\\[.2em]
\displaystyle\textrm{such that $\WhI\subseteq \nos\inv(Y)$}}
\ \right\}
\ \overset{\displaystyle\Psi}{\underset{\simeq}{\,\longleftarrow\,}}\ 
\left\{\ 
\substack{\displaystyle\textrm{Serre $\otimes$-ideals~$\cJ$ of~$\mL$}
\\[.2em]
\displaystyle\textrm{such that $\WhJ\subseteq \nos\inv(Y')$}}
\ \right\}
\]
given by $\Psi(\cJ)=(F_!)\inv(\cJ)$.
Moreover, if~$\cI=\Psi(\cJ)$ then $F$ induces an equivalence
\begin{equation}
\label{eq:main}%
\frac{\MK}{\Loc(\cI)}\isotoo\frac{\ML}{\Loc(\cJ)}
\end{equation}
between the~$\cI$-layer of~$\cK$ and the $\cJ$-layer of~$\cL$.
\end{thm}

\begin{proof}
The preimage of a localizing $\otimes$-ideal under a cocontinuous exact tensor-functor of tensor Grothendieck categories is again a localizing $\otimes$-ideal. Using~$F_!\colon \MK\to \ML$ and~$j^*\circ F_!\cong (F_Y)_!\circ i^*$ from~\Cref{lem:F!-restr}, we get a commutative diagram (of posets)
\begin{equation}
\label{eq:loc-square}%
\vcenter{
\xymatrix@C=5em@R=2em{
\left\{\ 
\substack{\displaystyle\textrm{Localizing $\otimes$-ideals of~$\MK$}
\\[.1em]
\displaystyle\textrm{containing~$\Bmin(Y)=\Ker(i^*)$}}
\ \right\}
& \left\{\ 
\substack{\displaystyle\textrm{Localizing $\otimes$-ideals of~$\ML$}
\\[.1em]
\displaystyle\textrm{containing~$\Bmin(Y')=\Ker(j^*)$}}
\ \right\}
\ar[l]_-{(F_!)\inv(-)}^-{=:\Psi}
\\
\left\{\ 
\textrm{Localizing $\otimes$-ideals of~}\MK_Y
\ \right\}
\ar[u]_-{\simeq}^-{(i^*)\inv(-)}
\ar@{<-}[r]^{((F_Y)_!)\inv(-)}_-{\simeq}
& 
\left\{\ 
\textrm{Localizing $\otimes$-ideals of~}\ML_{Y'}
\ \right\}
\ar[u]^-{\simeq}_-{(j^*)\inv(-)}
}}\kern-1em
\end{equation}
involving the big categories of (not necessarily finitely presented) modules.
The bottom map is a bijection since~$F_Y=F\restr{\cK_Y}\colon \cK_Y\isoto \cL_{Y'}$ is an equivalence. The vertical maps are bijections since $\MK_Y=(\MK)/\Bmin(Y)$ by~\Cref{lem:ker(i*)} and similarly on the~$\cL$-side.
Let us denote by~$\Psi$ the resulting bijection $\Psi(\cC)=(F_!)\inv(\cC)$ at the top of~\eqref{eq:loc-square}. 
(We shall show that it `restricts' to the $\Psi$ of the statement.)
Whenever~$\cB$ corresponds to~$\cC$ under this bijection, the functor $F_!\colon \MK\to \ML$ descends to an equivalence
\begin{equation}
\label{eq:equiv-Phi}%
\frac{\MK}{\cB}\isotoo \frac{\ML}{\cC}.
\end{equation}
Indeed, by the above discussion (and \Cref{lem:F!-restr} again) this functor is simply the equivalence~$(F_Y)_!$ descended to the corresponding quotients of~$\MK_Y$ and~$\ML_{Y'}$. Since this is a $\otimes$-functor (as~$F_!$ is) the $\otimes$-unit is finitely presented on one side of~\eqref{eq:equiv-Phi} if and only if it is on the other.
In other words, the $\Psi$ of~\eqref{eq:loc-square} restricts to a bijection on the \emph{finite type} localizing $\otimes$-ideals by \Cref{thm:finite-type}. This gives the bijection at the top of the following diagram:
\begin{equation*}
\label{eq:fin-loc-rectangle}%
\vcenter{
\xymatrix@C=4em@R=2em{
\left\{\ 
\substack{\displaystyle\textrm{Localizing $\otimes$-ideals of~$\MK$}
\\[.1em]
\displaystyle\textrm{of \emph{finite type} containing~$\Bmin(Y)$}}
\ \right\}
\ar[d]^-{(-)\fp=(-)\cap\mK}_-{\simeq}
& \left\{\ 
\substack{\displaystyle\textrm{Localizing $\otimes$-ideals of~$\ML$}
\\[.1em]
\displaystyle\textrm{of \emph{finite type} containing $\Bmin(Y')$}}
\ \right\}
\ar[l]_-{\Psi}^-{\simeq}
\ar[d]^-{(-)\fp=(-)\cap\mL}_-{\simeq}
\\
\left\{\ 
\substack{\displaystyle\textrm{Serre $\otimes$-ideals $\cI$ of~$\mK$}
\\[.2em]
\displaystyle\textrm{such that $\WhI\subseteq\nos\inv(Y)$}}
\ \right\}
& \left\{\ 
\substack{\displaystyle\textrm{Serre $\otimes$-ideals $\cJ$ of~$\mL$}
\\[.2em]
\displaystyle\textrm{such that $\WhJ\subseteq\nos\inv(Y')$}}
\ \right\}
\ar@{..>}[l]_-{}
}}
\end{equation*}
The vertical maps on both sides are well-defined bijections by~\Cref{thm:Bmin} and what it means to be of finite type. (Their inverses are given by~$\Loc(-)$.)
So we get the (dotted) bijection at the bottom and we claim it is the~$\Psi$ of the statement.
Indeed for every localizing $\otimes$-ideal~$\cC$ of~$\ML$, which is of finite type and contains~$\Bmin(Y')$, the following Serre $\otimes$-ideals of~$\mK$ coincide because $F_!$ preserves finitely presented objects:
\[
\Psi(\cC)\cap\mK=\big((F_!)\inv(\cC)\big)\cap\mK=\big((F_!)\inv(\cC\cap \mL)\big)\cap \mK.
\]
This gives the bijection~$\Psi$ announced in the statement (where we only used~$F_!$ on finitely presented modules).
Finally,~\eqref{eq:main} follows from~\eqref{eq:equiv-Phi} for~$\cB=\Loc(\cI)$ and~$\cC=\Loc(\cJ)$.
\end{proof}
\begin{rem}
One can trace through the above proof to show that~$\Psi\inv$ maps a Serre $\otimes$-ideal~$\cI$ to the finitely presented part of the localizing $\otimes$-ideal generated by~$F_!(\cI)\cup\Bmin(Y')$.
\end{rem}

Let us denote by~$\varphi\homol=\Spch(F)\colon \SpchL\to \SpchK$ the map induced by~$F$ on homological spectra, given by preimages~$\varphi\homol(\cG)=(F_!)\inv(\cG)$ under~$F_!\colon \mK\to \mL$. Even in the absence of the steel property, the bijection $\varphi\inv(Y)\isoto Y$ of \Cref{prop:Y-equivalence} extends to~$\Spch$.
\begin{cor}
\label{cor:Y-equiv-Spch}%
Let $F\colon \cK\to \cL$ be a tt-functor that is an equivalence over a Thomason subset~$Y\subseteq\SpcK$ and let $Y'\subseteq\SpcL$ be its preimage. Then $\varphi\homol\colon \SpchL\to \SpchK$ restricts to a bijection $\nos\inv(Y')\isoto \nos\inv(Y)$. Moreover, for every~$\cG\in\nos\inv(Y')$ the functor~$F$ induces an equivalence between the homological residue fields at~$\cG$ and at~$\cH=\varphi\homol(\cG)$
\[
\frac{\MK}{\Loc(\cH)}\isotoo \frac{\ML}{\Loc(\cG)}.
\]
\end{cor}
\begin{proof}
Apply \Cref{thm:main}. When~$\cI=\Psi(\cJ)$, the equivalence~\eqref{eq:main} on the  quotients shows that $\cI$ is maximal (a homological prime) if and only if~$\cJ$ is maximal. 
For~$\cH\in\SpchK$ we have $\Wh(\cH)=\{\cH\}$ and so the condition~$\Wh(\cH)\subseteq \nos\inv(Y)$ simply reads $\cH\in\nos\inv(Y)$. Similarly, $\Wh(\cG)\subseteq \nos\inv(Y')$ if and only if~$\cG\in\nos\inv(Y')$. 
Since $\Psi(\cG)=(F_!)\inv(\cG) = \phih(\cG)$, we see that the bijection~$\Psi$ restricts to the claimed bijection~$\varphi\homol\colon\nos\inv(Y')\isoto \nos\inv(Y)$.
\end{proof}

\begin{rem}
Since the bijection~$\Psi$ of \Cref{thm:main} preserves inclusions and restricts to~$\varphi\homol$ on maximal ideals, we see that whenever~$\cI=\Psi(\cJ)$ then~$\varphi\homol(\WhJ)=\WhI$, and thus
\begin{equation}
\label{eq:phih-continuous}
(\varphi\homol)\inv(\WhI)=\WhJ.    
\end{equation}
\end{rem}
\begin{rem}
\label{rem:main-with-equality}%
The equality \eqref{eq:phih-continuous} gives a specialization of \Cref{thm:main} when $Y$ (hence~$Y'$) is proconstructible.
Since~$(\phih)\inv(\nos\inv(Y)) = \nos\inv(Y')$, the bijection~$\Psi$ restricts to a bijection between those $\cI$ with $\WhI = \nos\inv(Y)$ and those $\cJ$ with $\WhJ = \nos\inv(Y')$.
\end{rem}

Here is another formulation of the above preservation of homological residue fields, that we shall need in \Cref{sec:regularity}. We keep notation as above.

\begin{cor}
\label{cor:iso-residue}%
Let $F\colon \cK\to \cL$ be a tt-functor that is an equivalence over~$Y$. Set $Y' = \varphi\inv(Y)$. Let $\cJ\subseteq\mL$ be a Serre $\otimes$-ideal such that $\WhJ\subseteq\nos\inv(Y')$ and~$\cI=(F_!)\inv(\cJ)$ its correspondent in~$\mK$. Then we have an isomorphism~$F^*(E_{\cJ})\simeq \EI$ in~$\MK$. In particular, if $\cG \in \nos \inv(Y')$ in~$\SpchL$ and $\cH \coloneq\varphi\homol(\cG)$ then $F^*(\kappa(\cG))\simeq\kappa(\cH)$.
\end{cor}
\begin{proof}
By \Cref{thm:main}, we have an equivalence of layers, fitting in the commutative diagram
\[
\xymatrix@R=1.5em{
\MK \ar[r]^-{F_!} \ar@{->>}[d]_-{\QI}
& \ML \ar@{->>}[d]^-{\QJ}
\\
\MK/\Loc(\cI) \ar[r]^-{\overline{F_!}}_-{\simeq}
& \ML/\Loc(\cJ).
}
\]
On right adjoints, we get $F^*\circ R_{\cJ}\cong \RI\circ (\overline{F_!})\inv$.
Applying this isomorphism to the object $\bar{E}_{\cJ}$ we get $F^*(E_{\cJ})\simeq \EI$ since the equivalence $\overline{F_!}$ identifies the injective hulls $\bar{E}_{\cI}$ and $\bar{E}_{\cJ}$.
\end{proof}

\begin{rem}
The bijection~$\Psi$ of \Cref{thm:main} preserves the Serre $\otimes$-ideals of the form~$\cI_{S}$ as in~\Cref{dfn:I_S}. Indeed, let $S\subseteq \cK$ be a $\otimes$-multiplicative class of objects in $\cK$. By~\Cref{cor:I_S} the homological spectrum of~$\cI_{S}$ is~$\nos\inv(Z(S))$ where $Z(S)=\cap_{a\in S}\supp(a)$ is the closed subset of $\Spc(\cK)$ corresponding to~$S$. Suppose that $Z(S)$ is a Thomason subset and that $F\colon \cK \to \cL$ is an equivalence over $Z(S)$. As $F$ is a tt-functor, the set $F(S)$ is a $\otimes$-multiplicative set of objects in $\cL$ and $Z(F(S)) = \varphi\inv(Z(S))$. We are claiming that
\[
\Psi(\cI_{F(S)})=(F_!)\inv(\cI_{F(S)}) = \cI_S. 
\]
Using \Cref{prop:I_S}, the inclusion $\cI_S \subseteq (F_!)\inv (\cI_{F(S)})$ is clear. For the reverse inclusion, suppose that $M \in \mK$ and $F_!(M) \in \cI_{F(S)}$. This means that $0=\yon(F(a))\otimes F_!(M)\cong F_!(\yon(a) \otimes M)$ for some~$a\in S$. We may assume that~$\supp(a) = Z(S)$ by replacing~$a$ with~$a_0\otimes a$, where~$a_0\in S$ and $\supp(a_0)=Z(S)$ as in \Cref{rem:I_S-Thomason}. Then $F_!(\yon(a) \otimes M) =0$ in~$\ML$ forces $i^*(\yon(a) \otimes M) =0$ in~$\Modd{\cK_{Z(S)}}$ by \Cref{lem:F!-restr} applied to~$Y=Z(S)$. Since $\yon(a) \otimes M$ is $\Ker(i^{\ast})$-local (\Cref{lem:ker(i*)}), we get $\yon(a) \otimes M =0$. In other words, $M\in\cI_S$ as claimed.
\end{rem}

To conclude this section, we apply our results to completion.
\begin{rec}
\label{rec:completion}%
Suppose that $\cK$ comes from an ambient `big' tt-category~$\cT$, as in~\eqref{eq:K-T}. The \emph{Bousfield completion} of $\cT$ along a Thomason $Y \subseteq \SpcK$ is the subcategory $\cT_Y^{\perp\perp}=\big((\cK_Y)^{\perp}\big)^{\perp}$ of~$Y$-complete objects in~$\cT$.
This should not be mistaken for \emph{the} completion~$\hT$ as a `big' tt-category. See~\cite[Example~2.13 and Remark~5.4]{BalmerSanders26}.
Only the dualizables of~$\hT$ should be expected to be found in~$\cT_Y^{\perp \perp}$ (and thus ultimately in~$\cT$!) and one \emph{defines}
\[
\hT\dual\coloneq(\cT_Y^{\perp \perp})\dual.
\]
The full~$\hT$ is the Ind-completion of~$\hT\dual$.
Denote by $\eY\in\cT$ the left idempotent of \cite{BalmerFavi11}.
The localization $[\eY, -] \colon \cT \to \cT_Y^{\perp \perp}$ is a tensor functor and its restriction to dualizable objects is called the \emph{$Y$-completion functor} $\Compl_Y\colon \cT\dual\to \hT\dual$.
By \cite[Theorem 4.6]{BalmerSanders25}, it is an equivalence over~$Y$. In light of this, the following, announced as \Cref{thm:main-intro}, is immediate.
\end{rec}

\begin{cor}
\label{thm:residue-completion}
Keep notation as in~\Cref{rec:completion}. Then the $Y$\!-completion functor $\Compl_Y\colon \cT\dual\to \hT\dual$ induces an inclusion-preserving bijection between the Serre $\otimes$-ideals $\cI\subseteq \mT$ with $\WhI \subseteq \nos\inv(Y)$ and the Serre $\otimes$-ideals~$\hat\cI\subseteq \mhT$ with $\Wh(\hat{\cI}) \subseteq \nos\inv(\varphi\inv(Y))$, such that we have an equivalence of the corresponding layers
    \begin{equation*}
    \frac{\MT}{\Loc(\cI)} \isotoo \frac{\MhT}{\Loc(\hat\cI)}.
    \end{equation*}
In particular, if $\cI\in \nos\inv(Y)$ is maximal, then $\hat\cI \in \nos\inv(\varphi\inv(Y))$ is maximal as well and the above is an equivalence of homological residue fields.
\qed\end{cor}

\section{Residual regularity}
\label{sec:regularity}%

Let~$\cK$ and~$\cL$ be tt-categories as in our \Cref{hyp:K-T} and~$F\colon \cK\to \cL$ be a tt-functor that is an equivalence over a Thomason subset~$Y$ (\Cref{dfn:Y-equivalence}). We want to show that~$F$ preserves and detects residual regularity on~$Y$. As usual, we write~$\varphi\coloneq\Spc(F)$ and~$\varphi\homol\coloneq\Spch(F)$.

\begin{rec}
\label{rec:stalkcats}
For any prime $\otimes$-ideal $\cP$ of $\cK$, we define the \emph{local category} at~$\cP$ to be $\cK_{\cP} \coloneq (\cK/\cP)^{\natural}$, which is again an essentially small, rigid and idempotent complete tt-category. By \cite[Proposition 2.9]{Balmer20a} the functor $q_! \colon \mK \to \mK_{\cP}$ induced by $q \colon \cK \to \cK_{\cP}$ is a Gabriel quotient with kernel~$\ideal{ \yon(\cP)}$. It follows that $\cG \mapsto (q_!)\inv(\cG)$ defines an inclusion-preserving one-to-one correspondence $\Spch(\cK_\cP)\isoto \Wh(\ideal{\yon(\cP)})$ between homological primes $\cG$ of $\cK_{\cP}$ and the homological primes $\cH$ of $\cK$ which contain the Serre $\otimes$-ideal $\ideal{ \yon(\cP) }$. Moreover, if $\cG$ corresponds to $\cH$, then the respective homological residue fields are canonically equivalent in such a way that the following diagram is commutative; see \cite[Corollary 3.6]{Balmer20a}: 
\[\begin{tikzcd}[row sep=1.5em]
	\cK & \MK & {\MK/\Loc(\cH)} \\
	{\cK_{\cP}} & {\Modd{\cK_{\cP}}} & {\Modd{\cK_{\cP}}/\Loc(\cG)}.
	\arrow["h", hook, from=1-1, to=1-2]
	\arrow["q"', two heads, from=1-1, to=2-1]
	\arrow["{Q_{\cH}}", two heads, from=1-2, to=1-3]
	\arrow["{q_!}", two heads, from=1-2, to=2-2]
	\arrow[from=1-3, to=2-3]
	\arrow["\sim"{marking, allow upside down}, shift left=2, draw=none, from=1-3, to=2-3]
	\arrow["h"', hook, from=2-1, to=2-2]
	\arrow["{Q_{\cG}}"', two heads, from=2-2, to=2-3]
\end{tikzcd}\]
Since the injective hulls $\bar{E}_{\cH}$ and $\bar{E}_{\cG}$ are identified under the right-hand equivalence, we have $\kappa(\cH)\simeq q^*(\kappa(\cG))$ in $\MK$ and $q_! (\kappa(\cH)) \simeq \kappa(\cG)$ in $\Modd{\cK_{\cP}}$. In particular, $\kappa(\cH)$ is a $\Ker(q_!)$-local object and it makes sense to speak of $\kappa(\cH)$ as an object in $\Modd{\cK_{\cP}}$. 
\end{rec}

The following definition is essentially \cite[Definition~2.1]{VanRooy26}. 
\begin{dfn}\label{dfn:resreg}
The tt-category $\cK$ is called \emph{residually regular at the point ${\cP \in \SpcK}$} if there exists a homological prime $\cH$ in $\nos\inv(\{\cP\})$ and an object~$a$ in the local category~$\cK_{\cP}=(\cK/\cP)^\natural$ whose image under Yoneda $\yon\colon \cK_\cP\hook\Modd{\cK_\cP}$ is isomorphic to~$\kappa(\cH)$.
\end{dfn}
\begin{rem}
In~\cite{VanRooy26}, the tt-category $\cK$ is assumed to be the small part of a `big' tt-category as in~\eqref{eq:K-T}. This allows us to produce~$\kappa(\cH)$ in~$\cT$ by Brown Representability (see~\cite[Corollary~1.9]{Krause00}). Since an object~$X\in\cT$ becomes isomorphic to a compact~$a\in\cT\compact$ in~$\Modd{\cT\compact}$ if and only $X\simeq a$ in~$\cT$, \Cref{dfn:resreg} coincides with~\cite[Definition~2.1]{VanRooy26} in this case.
\end{rem}

\begin{rem}\label{rem:resreg-stalklocal}
From \Cref{rec:stalkcats}, it follows that $\cK$ is residually regular at the point ${\cP \in \SpcK}$ if and only the local category~$\cK_{\cP}$ is residually regular at its unique closed point~$(0)$. 
\end{rem}

\begin{prop}
\label{prop:residueiscompact}%
Suppose that $\cK$ is a local tt-category.
Then~$\cK$ is residually regular at the closed point~$(0)$ if and only if there exists a homological prime in~$\cH\in\nos \inv \{(0) \}$ such that $\kappa(\cH)$ is finitely presented in $\MK$.
In that case, the unique object~$a\in\cK$ such that $\yon(a)\simeq \kappa(\cH)$ has support exactly the closed point, $\supp(a) = \{ (0)\}$, and~$\nos\inv(\{(0)\})=\{\cH\}$.
\end{prop}
\begin{proof}
If~$\kappa(\cH)$ is finitely presented then it is an injective object in~$\mK$. Hence it must come from a unique~$a \in \cK$, that is, $\yon(a) \simeq\kappa(\cH)$.
For every~$\cG\in\SpchK$ the object $\kappa(\cG)\otimes \yon(a)\simeq\kappa(\cG)\otimes\kappa(\cH)$ is zero if and only if~$\cG\neq \cH$ (see \cite[Proposition~5.3]{Balmer20a}). Hence $\nos\inv(\supp(a))=$ $\SET{\cG\in\SpchK}{\kappa(\cG)\otimes \yon(a)\neq 0}=\{\cH\}$. By surjectivity of~$\nos$, we get $\supp(a)=\nos(\nos\inv(\supp(a)))=\nos(\{\cH\})=\{(0)\}$ and~$\nos\inv(\{(0)\})=\nos\inv(\supp(a))=\{\cH\}$.
\end{proof}

\begin{lem}
\label{lem:resreg-in-K_Y}
    Suppose that $\cK$ is a local tt-category and let $Y\subseteq \SpcK$ be a nonempty Thomason subset and $\cH$ a homological prime in $\nos \inv \{(0) \}$. Then, the residue object $\kappa(\cH)$ is finitely presented in $\MK$ if and only if $i^{\ast}(\kappa(\cH))$ is finitely presented in~$\MK_Y$.
\end{lem}
\begin{proof}
    First, suppose that $\kappa(\cH)$ is finitely presented in~$\MK$. By \Cref{prop:residueiscompact}, there exists an object $a \in \cK$ with $\supp(a) = \{ (0)\}$ and such that $\kappa(\cH) \simeq \yon(a)$. Since $Y$ is a nonempty Thomason subset, it must contain the closed point $(0)$ and therefore $a \in \cK_Y$. It follows that $i^{\ast} (\kappa(\cH) )\cong i^{\ast}(\yon(a)) = \Hom_{\cK_Y}(-, a)$ is finitely presented in~$\MK_Y$. 

Conversely, suppose that $i^{\ast}(\kappa(\cH))$ is finitely presented in~$\MK_Y$. 
Since $i^{\ast}$ has an exact left adjoint $i_!$ (see \Cref{rec:F!-F*}), it preserves injective objects. 
So, $i^{\ast}(\kappa(\cH))$ is in fact an injective object in $\mK_Y$, therefore $i^*(\kappa(\cH))\simeq\yon(a)$ for some $a\in\cK_Y$. 
Note that $\kappa(\cH)$ is a $\Ker(i^{\ast})$-local object since it is $\Loc(\cH)$-local, and we have $\Ker(i^{\ast}) \subseteq \Loc(\cH)$ by \Cref{thm:Bmin} (because $\nos(\cH) =(0) \in Y$). Therefore $\kappa(\cH)\cong i_!i^*(\kappa(\cH))\simeq i_!(\yon(a))\cong\yon(i(a))\in\mK$.
\end{proof}

\begin{thm}
\label{thm:regularity}%
    Let $F\colon \cK\to \cL$ be a tt-functor. Let $Y\subseteq \SpcK$ be a nonempty Thomason subset such that $F$ is an equivalence over~$Y$. Let $\cQ \in \varphi\inv(Y)$ and set $\cP = \varphi(\cQ)$. Then $\cL$ is residually regular at $\cQ$ if and only if $\cK$ is residually regular at~$\cP$. 
\end{thm}
\begin{proof}
Since $\cP=F\inv(\cQ)$ we get a conservative localization $\bar{F}\colon \cK_{\cP}\to \cL_{\cQ}$ as in \Cref{prop:excision} applied to~$U=\Supp(\cP)\compl$ and~$V=\Supp(\cQ)\compl$. By that proposition, we can assume that $\cK$ and $\cL$ are local tt-categories and $\cP$ and $\cQ$ are their respective closed points. Since residual regularity is a local property (see \Cref{rem:resreg-stalklocal}), this suffices to prove the statement.

\Cref{cor:Y-equiv-Spch} implies that $\varphi\homol$ restricts to a bijection $\nos\inv(\{\cQ\})\isoto \nos\inv(\{\cP\})$. Consider a matching pair of points $\cG\in\nos\inv(\{\cQ\})$ and~$\cH=\varphi\homol(\cG)$. Writing as before $i\colon \cK_Y\hook \cK$ and~$j\colon \cL_{Y'}\hook\cL$ for the inclusions we have $j\circ F_Y=F\circ i$ and therefore $(F_Y)^*j^*(\kappa(\cG))\cong i^*F^*(\kappa(\cG))\simeq i^*(\kappa(\cH))$, where the last isomorphism uses~\Cref{cor:iso-residue}.
Since $(F_Y)^*\colon$ $\Modd{\cL_{Y'}}\isoto \Modd{\cK_{Y}}$ is an equivalence, it follows that $j^*(\kappa(\cG))$ is finitely presented if and only if $i^*(\kappa(\cH))$ is.
We can now translate these properties into the residual regularity of~$\cK$ at~$\cP=(0)$, respectively of~$\cL$ at~$\cQ=(0)$ thanks to~\Cref{lem:resreg-in-K_Y} and~\Cref{prop:residueiscompact}.
\end{proof}

\Cref{thm:regularity} applies in particular to the completion as in \Cref{rec:completion}. 
\begin{cor}
    Let $\cT$ be a rigidly-compactly generated tt-category, and let $Y \subseteq \SpcT$ be a nonempty Thomason subset. Let $\hat{\cP} \in \varphi\inv(Y)$ and set $\cP = \varphi(\hat{\cP})$. Then $\hT\dual$ is residually regular at $\hat{\cP}$ if and only if $\cT\dual$ is residually regular at~$\cP$. 
\qed
\end{cor}

\section{Example: Commutative rings}
\label{sec:commutative-rings}%

Let $R$ be a commutative ring. In this section, we consider layers of $\cK \coloneq \Dperf(R)$. In particular, we study the interplay of ideals $I \subseteq R$ and Serre $\otimes$-ideals $\cI \subseteq \mK$.

\begin{rec}
The tt-category $\Dperf(R)$ has the steel property~$\nos\colon \SpchK\isoto \SpcK$ by \cite[Corollary 5.11]{Balmer20a}. Its tt-spectrum identifies with the Zariski spectrum of the ring~$R$ via the comparison map $\comp\colon\SpcK \isoto \Spec(R)$ defined by~$\comp(\cP)=\SET{r\in R}{\cone(r)\notin\cP}$. See~\cite{Balmer10b}.
We sometimes drop the homeomorphisms~$\nos$ and~$\comp$ for readability.
The comparison map satisfies $\comp\inv(V(r))=\supp(\cone(r))$ for every $r\in R$ and consequently
\begin{equation}
\label{eq:supp-kos}%
\comp\inv(V(r_1,\ldots,r_n))=\supp(\kos(\underline{r}))
\qquadtext{where}
\kos(\underline{r}) = \otimes_{i=1}^n \cone(r_i)
\end{equation}
is the Koszul complex associated to the sequence~$\underline{r}=(r_1,\ldots,r_n)$.
\end{rec}

\begin{cons}
\label{cons:ideal-Serre-ideal-corr}
Given an ideal $I \subseteq R$, let $S_I = \SET{\kos(\underline{r})}{\underline{r}\textrm{ in }I}$ be the $\otimes$-multiplicative subset of~$\cK=\Dperf(R)$ generated by~$\SET{\cone(r)}{r\in I}$. Denote the Serre $\otimes$-ideal $\cI_{S_I}$ of \Cref{dfn:I_S} by $\sigma(I)$. The reader can verify that $\ker(\yon(\coev_{\cone(r)})) \cong \im(\yon(r))$ for every element~$r \in R$, so the Serre $\otimes$-ideal~$\sigma(I)$ is generated by $\SET{\im(\yon(r))}{r \in I}$ in~$\mK$. By \Cref{prop:I_S}, we also have the following description
\[\sigma(I) = \SET{M \in \mK}{\yon(\kos(\underline{r})) \otimes M = 0 \textrm{ for some } \underline{r}=(r_1, \ldots, r_n)\textrm{ in }I}.\]

In the other direction, given a Serre $\otimes$-ideal $\cI \subseteq \mK$, define an ideal of~$R$
\[
\tau(\cI) \coloneq \SET{r\in R}{\im(\yon(r))\in\cI}=\ker\big(\End_\cA(\unit)\xto{\,\QI\,} \End_{\AI}(\unit)\big)
\]
using that $R\cong\End_{\cK}(\unit)\cong\End_{\cA}(\unit)$ and the exact $\otimes$-functor~$\QI\colon \cA=\MK\onto \AI$.
\end{cons}

\begin{prop}
The functions~$\sigma$ and~$\tau$ of \Cref{cons:ideal-Serre-ideal-corr} are compatible with vanishing loci in the following sense.
\begin{enumerate}[\rm(a)]
\item
\label{it:spc-sigma)}%
For every ideal~$I\subseteq R$, we have $\comp\nos(\Wh(\sigma(I)))=V(I)$ in~$\Spec(R)$. Omitting the homeomorphisms for readability, this amounts to
\[
\W(\sigma(I))=V(I).
\]
In particular, the $\sigma(I)$-layer is closed.
\smallbreak
\item
\label{it:V(tau)}%
For every Serre $\otimes$-ideal~$\cI\subseteq \mK$, the closure of~$\comp\nos(\Wh(\cI))$ in~$\Spec(R)$ is $V(\tau(\cI))$. Omitting the homeomorphisms again, this reads
\[
\overline{\WI}=V(\tau(\cI)).
\]
In particular, we have $\W(\cI)\subseteq V(\tau(\cI))$, with equality if and only if the~$\cI$-layer is closed.
\end{enumerate}
\end{prop}
\begin{proof}
\eqref{it:spc-sigma)}:
We use the characterization $\sigma(I) = \cI_S$, where $S = \SET{\kos(\underline{r})}{\underline{r} \text{ in } I}$. By \Cref{cor:I_S}, the spectrum~$\W(\sigma(I))=\W(\cI_S)$ is
\[
\bigcap_{a\in S}\supp(a)=
\!\bigcap_{\underline{r}\textrm{ in }I}\!\supp(\kos(\underline{r}))\overset{\textrm{(\ref{eq:supp-kos})}}{=}
\!\bigcap_{\underline{r}\textrm{ in }I}\!\comp\inv(V(\underline{r}))=
\comp\inv\!\bigcap_{\underline{r}\textrm{ in }I}V(\underline{r})=
\comp\inv V(I).
\]
\eqref{it:V(tau)}: 
\textit{Claim: For~$r\in R$, we have $\comp\W(\cI)\subseteq V(r)$ if and only if~$r\in\sqrt{\tau(\cI)}$.}
Suppose first that $\comp\W(\cI)\subseteq V(r)$. This means that for every~$\cH\supseteq \cI$, we have $r\in\comp\nos(\cH)$, that is, $\overline{\cone(r)}\neq0$ in~$\AH$. Since $r\otimes \cone(r)=0$, we must have~$\overline{r}=0$ in~$\AH$ for every $\cH\supseteq\cI$. By \Cref{rec:nilpotence}, this forces~$r\colon \unit\to \unit$ to become $\otimes$-nilpotent in~$\AI$, meaning~$r\in\sqrt{\tau(\cI)}$. Conversely, if~$r^n\in\tau(\cI)$ for~$n\ge 1$, then $r^n=0$ in~$\AI$ and a fortiori in~$\AH$ for every~$\cH\in \Wh(\cI)$. This means that $\cone(r^n)$ goes to a non-zero object in~$\AH$ and thus~$r^n \in \comp \nos(\cH)$, forcing $r \in \comp\nos(\cH)$ by primality.
This proves the Claim. Since~$\{V(r)\}_{r\in R}$ forms a closed basis of the Zariski topology, the Claim implies that the closure of~$\comp\W(\cI)$ is the intersection of all the $V(r)$ for~$r\in\sqrt{\tau(\cI)}$. This intersection is~$V(\sqrt{\tau(\cI)})=V(\tau(\cI))$.
\end{proof}

The next proposition shows that the collection of closed layers of $\Dperf(R)$ carries at least as much information as the collection of quotients of~$R$.
\begin{prop}
\label{prop:tau-retraction-of-sigma}%
For every ideal $I \subseteq R$, we have $\tau(\sigma(I)) = I$. So we have an injection \[\sigma \colon \{\textrm{ideals of } R\} \to \{\textrm{closed layers of } \Dperf(R)\}.\]
\end{prop}
\begin{proof}
We have $\tau(\sigma(I)) = \SET{x \in R}{\im(\yon(x)) \in \sigma(I)}.$ Certainly, $I \subseteq \tau(\sigma(I))$, since $\cone(x) \otimes x = 0$ for any $x \in I$. For the other inclusion, if~$x \in R$ satisfies $\yon(\kos(\underline{r})) \otimes \im(\yon(x)) = 0$ for $\underline{r}=(r_1, \ldots, r_n)$ in~$I$, then $\kos(\underline{r}) \otimes x = 0$. So multiplication by $x$ is zero on the homology of $\kos(\underline{r})$, including $H_0(\kos(\underline{r})) = R / \ideal{r_1,\ldots,r_n}$. This implies $x \in \ideal{r_1,\ldots,r_n} \subseteq I$.
\end{proof}

\begin{rem}
In general, one cannot expect $\sigma$ to be a bijection. For example, let $R = \bbZ$. The only ideal~$I$ with $V(I) = \Spec R$ is $I=(0)$. However, there are many Serre $\otimes$-ideals~$\cI$ other than $\sigma((0)) = (0)$ for which $\Wh(\cI) = \SpchK$. For example, let $c$ be the complex $\bbZ \xto{4} \bbZ$, and consider the map $f = (c \xto{2} c)$. Then $f$ vanishes in every residue field, hence $\cI \coloneq \ideal{\im(\yon(f))}$ is contained in every $\cH \in \SpchK$. So~$\W(\cI)=\Spec(\bbZ)$ and~$\tau(\cI)$ must be~$0$, but~$\cI\neq 0$.
\end{rem}

\begin{rem}
Given an ideal $I \subseteq R$, there is another Serre $\otimes$-ideal of $\mK$ associated to~$I$, namely the kernel of the functor $F_! \colon \mK \to \modd \Dperf(R / I)$ induced by the ring map $R \onto R / I$. This kernel contains $\sigma(I)$ but can be larger in general. For example, let $k$ be a field of characteristic $0$, let $R = k[t] / (t^2)$, and let $f$ be the following morphism of complexes:
\[\begin{tikzcd}[row sep=1em]
    \cdots & 0 & R & R & 0 & \cdots \\
    \cdots & 0 & 0 & R & 0 & \cdots
    \arrow[from=1-1, to=1-2]
    \arrow[from=1-2, to=1-3]
    \arrow["t", from=1-3, to=1-4]
    \arrow[from=1-4, to=1-5]
    \arrow[from=1-5, to=1-6]
    \arrow[from=2-1, to=2-2]
    \arrow[from=2-2, to=2-3]
    \arrow[from=2-3, to=2-4]
    \arrow[from=2-4, to=2-5]
    \arrow[from=2-5, to=2-6]
    \arrow[from=1-2, to=2-2]
    \arrow[from=1-3, to=2-3]
    \arrow[from=1-4, to=2-4]
    \arrow["\ t", from=1-4, to=2-4]
    \arrow[from=1-5, to=2-5]
\end{tikzcd}\]
Letting $I = (t)$, we see that $\im(\yon(f)) \in \Ker(F_!)$. However, $f \otimes \cone(t) \neq 0$. In fact, one can verify directly that no nonzero compact object $c \in \Dperf(R)$ annihilates~$f$ (using that~$\textrm{char}(k)=0$). See \cite{AmodeoMoschetti15} for details on the category $\Dperf(k[t] / (t^2))$.
\end{rem}



\providecommand{\bysame}{\leavevmode\hbox to3em{\hrulefill}\thinspace}
\providecommand{\MR}{\relax\ifhmode\unskip\space\fi MR }
\providecommand{\MRhref}[2]{%
  \href{http://www.ams.org/mathscinet-getitem?mr=#1}{#2}
}
\providecommand{\href}[2]{#2}


\end{document}